\documentclass[11pt]{amsart}
\usepackage[T1]{fontenc}
\usepackage[utf8]{inputenc}
\usepackage{lmodern}
\usepackage{microtype}
\usepackage{eucal}
\ifdefined\pdfgentounicode
\pdfglyphtounicode{minus}{002D}
\fi
\usepackage{amsmath,amssymb,amsthm}
\usepackage[a4paper,margin=1in]{geometry}
\usepackage{enumitem}
\usepackage{booktabs}
\usepackage{tikz}
\usepackage{xcolor}
\usepackage[normalem]{ulem}
\usetikzlibrary{decorations.pathreplacing}
\usepackage[hidelinks,pdftitle={First Order Geometry of Equilateral Small Polygons},pdfauthor={José Ayala Hoffmann, Tomás Ayala Coppola},pdfkeywords={Small polygons, equilateral polygons, maximal perimeter, first order optimality conditions, interval arithmetic, computer-assisted proof}]{hyperref}
\numberwithin{equation}{section}
\newtheorem{theorem}{Theorem}[section]
\newtheorem{lemma}[theorem]{Lemma}
\newtheorem{proposition}[theorem]{Proposition}
\theoremstyle{definition}
\newtheorem{definition}[theorem]{Definition}
\newtheorem{remark}[theorem]{Remark}
\theoremstyle{plain}

\newcommand{\R}{\mathbb R}
\newcommand{\Q}{\mathbb Q}
\newcommand{\Kcal}{\mathcal K}
\newcommand{\Ical}{\mathcal I}

\DeclareMathOperator{\rank}{rank}
\title{First Order Geometry of Equilateral Small Polygons}
\author{Jos\'e Ayala Hoffmann}
\address{Universidad de Tarapac\'a, Iquique, Chile}
\email{jayalhoff@gmail.com}

\author{Tom\'as Ayala Coppola}
\address{Universidad de Tarapac\'a, Iquique, Chile}
\email{tomasayalacoppola07@gmail.com}
\date{}
\subjclass[2020]{Primary 52A40; Secondary 52A10, 52C25}
\keywords{Small polygons, equilateral polygons, maximal perimeter,
first order optimality conditions, infinitesimal rigidity,
interval arithmetic, computer assisted proofs}
\begin{document}

\begin{abstract}
\baselineskip=20 true pt
\baselineskip=1.1\normalbaselineskip
We develop a first order calculus based on the geometry of realised diameter graphs for perimeter maximisation among convex equilateral small polygons. Polygon side equalities and active diameter constraints define a linearised feasible cone and an associated equality space. 
We prove that an equilateral small polygon is a
strict local maximiser of perimeter modulo Euclidean motions
if a multiplier balance holds with strictly positive active
diameter multipliers and the active equality space consists
only of infinitesimal Euclidean motions. The proof uses a local Euclidean slice and a geometric limiting argument and requires neither independence of the active constraint gradients nor second order information.

As an application, we combine our criterion with validated interval arithmetic to establish strict local maximality for Bingane and Audet's polygons $A_{16}$, $A_{32}$, $A_{64}$ and for reconstructions of their polygons $Z_{32}$ and $Z_{64}$. For each $n\in\{32,64\}$, the perimeter of $Z_n$ is strictly larger than that of $A_n$. Thus, strict local maximisers with different perimeters coexist in the space of convex equilateral small polygons with the same number of sides.
\end{abstract}
\maketitle

\baselineskip=20 true pt
\baselineskip=1.1\normalbaselineskip
\section{Introduction}

What is the largest possible perimeter of a convex $n$-gon whose
vertices are pairwise at distance at most one? This is known as the
isodiametric perimeter problem for small polygons. In 1922,
Reinhardt~\cite{reinhardt} proved that every convex small $n$-gon
$P$ satisfies
\[
L(P)\leq L_n=2n\sin\frac{\pi}{2n},
\]
where $L(P)$ denotes its perimeter. For $n$ odd, Reinhardt~\cite{reinhardt} showed that the regular
$n$-gon, whose sides and interior angles are all equal, is a
global maximiser. For $n$ even with an odd divisor, the results
of Reinhardt~\cite{reinhardt}, Vincze~\cite{vincze} and
Datta~\cite{datta} establish that every global maximiser is
equilateral but not regular. In contrast, Reinhardt~\cite{reinhardt}
proved that when $n$ is a power of two, the bound $L_n$ is not
attained. For $n=4$ and $n=8$, the maxima were determined by Tamvakis~\cite{tamvakis} and Audet, Hansen and Messine~\cite{audethansenmessine}, respectively. More recent developments include the zonogon based optimisation approach of Mulansky and Potschka~\cite{mulanskypotschka} and a preprint by Guo and Luo~\cite{guoluo2026}
proposing computer assisted proofs for the unrestricted
convex small polygon problem at $n=16$, $32$ and $64$.

The equilateral restriction leads to a distinct extremal problem. Among powers of two, the largest case for which the exact maximum
is known is $n=8$, solved by Audet, Hansen, Messine and
Perron~\cite{audethansenmessineperron} in 2004. For $n=2^k$ with $k\geq4$, the exact maximum
remains unknown. Mossinghoff~\cite{mossinghoff2008}
constructed convex equilateral small polygons for $n=2^k$, $k\geq4$,
whose perimeters approach Reinhardt's bound within $O(1/n^4)$.
Bingane and Audet~\cite{binganeaudet2022} denote this family by $M_n$
and construct an improved family $A_n$, satisfying
\[
L_n-L(A_n)=\frac{\pi^4}{n^4}+O\Bigl(\frac1{n^5}\Bigr),
\qquad
L(A_n)-L(M_n)=\frac{2\pi^4}{n^4}+O\Bigl(\frac1{n^5}\Bigr).
\]
Their constructions, together with OPTIGON ~\cite{optigon},
have provided a valuable basis for our work.

Our aim is to give sufficient geometric conditions for a
convex equilateral small polygon to be a strict local
maximiser of perimeter. Consider a strictly convex equilateral small
polygon and join every pair of vertices at unit distance.
The resulting realised diameter graph records the active diameter
constraints. Assign velocities to the vertices and allow the common
side length to vary at rate $\eta$. The side equalities and the
requirement that active diameters cannot increase define a
linearised feasible cone, and inside it the directions that preserve
every active diameter to first order. 
This treatment of active constraints is closely related to the
local calculus of realised contact geometry developed for the
minimum perimeter problem for disk packings~\cite{perimeterpaper}.

We prove that a strictly convex equilateral small polygon is a
strict local maximiser of perimeter modulo Euclidean motions
if a multiplier balance holds with strictly positive active
diameter multipliers and the active equality space consists
only of infinitesimal Euclidean motions. The proof is geometric: we remove translations and rotations by
placing nearby competitors in a local Euclidean slice. A sequence
of distinct competitors with at least equal perimeter would yield
a nonzero limiting direction after normalisation. The positive
multiplier balance forces this direction to preserve every active
constraint, while rigidity forces it to be an infinitesimal Euclidean
motion, a contradiction. Thus first order information at a single
configuration establishes strict local maximality modulo Euclidean
motions, without requiring a smooth feasible set or second order
information.

An important feature of our method is that the active constraint
gradients need not be independent, and the multipliers need not
be unique. The active diameter inequalities are analogous to
cable constraints in tensegrity theory~\cite{rothwhiteley,connelly}.

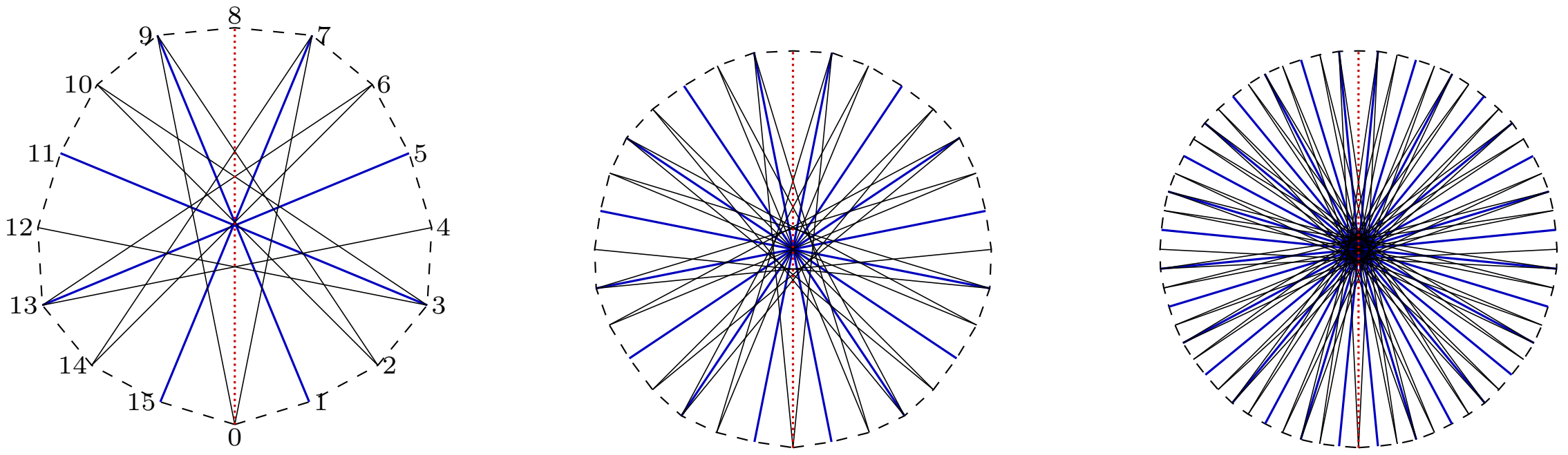
\begin{figure}[t]
\centering
\begin{minipage}[b]{0.32\textwidth}\centering
\begin{tikzpicture}[scale=3.8,
   outline/.style={dashed,line width=0.4pt},
   axis/.style={red!80!black,densely dotted,line width=0.6pt},
   path/.style={black,line width=0.3pt},
   pend/.style={blue!75!black,line width=0.6pt},
   lbl/.style={font=\tiny,inner sep=1pt}]
\coordinate (p0) at (0.00000,0.00000);
\coordinate (p1) at (0.18752,0.05686);
\coordinate (p2) at (0.36036,0.14920);
\coordinate (p3) at (0.48473,0.30063);
\coordinate (p4) at (0.49595,0.49626);
\coordinate (p5) at (0.43897,0.68375);
\coordinate (p6) at (0.34652,0.85653);
\coordinate (p7) at (0.19501,0.98080);
\coordinate (p8) at (0.00000,1.00000);
\coordinate (p9) at (-0.19501,0.98080);
\coordinate (p10) at (-0.34652,0.85653);
\coordinate (p11) at (-0.43897,0.68375);
\coordinate (p12) at (-0.49595,0.49626);
\coordinate (p13) at (-0.48473,0.30063);
\coordinate (p14) at (-0.36036,0.14920);
\coordinate (p15) at (-0.18752,0.05686);
\draw[outline] (p0) -- (p1) -- (p2) -- (p3) -- (p4) -- (p5) -- (p6) -- (p7) -- (p8) -- (p9) -- (p10) -- (p11) -- (p12) -- (p13) -- (p14) -- (p15) -- cycle;
\draw[pend] (p1) -- (p9);
\draw[pend] (p3) -- (p11);
\draw[pend] (p5) -- (p13);
\draw[pend] (p7) -- (p15);
\draw[path] (p0) -- (p7);
\draw[path] (p0) -- (p9);
\draw[path] (p2) -- (p9);
\draw[path] (p2) -- (p10);
\draw[path] (p3) -- (p10);
\draw[path] (p3) -- (p12);
\draw[path] (p4) -- (p13);
\draw[path] (p6) -- (p13);
\draw[path] (p6) -- (p14);
\draw[path] (p7) -- (p14);
\draw[axis] (p0) -- (p8);
\node[lbl,below] at (p0) {0};
\node[lbl,right] at (p1) {1};
\node[lbl,right] at (p2) {2};
\node[lbl,right] at (p3) {3};
\node[lbl,right] at (p4) {4};
\node[lbl,right] at (p5) {5};
\node[lbl,right] at (p6) {6};
\node[lbl,right] at (p7) {7};
\node[lbl,above] at (p8) {8};
\node[lbl,left] at (p9) {9};
\node[lbl,left] at (p10) {10};
\node[lbl,left] at (p11) {11};
\node[lbl,left] at (p12) {12};
\node[lbl,left] at (p13) {13};
\node[lbl,left] at (p14) {14};
\node[lbl,left] at (p15) {15};
\end{tikzpicture}

\small (a) \(A_{16}\)
\end{minipage}\hfill
\begin{minipage}[b]{0.32\textwidth}\centering
\begin{tikzpicture}[scale=3.8,
   outline/.style={dashed,line width=0.4pt},
   axis/.style={red!80!black,densely dotted,line width=0.6pt},
   path/.style={black,line width=0.3pt},
   pend/.style={blue!75!black,line width=0.6pt},
   lbl/.style={font=\tiny,inner sep=1pt}]
\coordinate (b0) at (0.00000,0.00000);
\coordinate (b1) at (0.09707,0.01440);
\coordinate (b2) at (0.19226,0.03824);
\coordinate (b3) at (0.28097,0.08020);
\coordinate (b4) at (0.35369,0.14610);
\coordinate (b5) at (0.41215,0.22492);
\coordinate (b6) at (0.46260,0.30909);
\coordinate (b7) at (0.49566,0.40148);
\coordinate (b8) at (0.49952,0.49954);
\coordinate (b9) at (0.48512,0.59661);
\coordinate (b10) at (0.46127,0.69180);
\coordinate (b11) at (0.41931,0.78051);
\coordinate (b12) at (0.35340,0.85322);
\coordinate (b13) at (0.27458,0.91168);
\coordinate (b14) at (0.19041,0.96213);
\coordinate (b15) at (0.09801,0.99518);
\coordinate (b16) at (0.00000,1.00000);
\coordinate (b17) at (-0.09801,0.99518);
\coordinate (b18) at (-0.19041,0.96213);
\coordinate (b19) at (-0.27458,0.91168);
\coordinate (b20) at (-0.35340,0.85322);
\coordinate (b21) at (-0.41931,0.78051);
\coordinate (b22) at (-0.46127,0.69180);
\coordinate (b23) at (-0.48512,0.59661);
\coordinate (b24) at (-0.49952,0.49954);
\coordinate (b25) at (-0.49566,0.40148);
\coordinate (b26) at (-0.46260,0.30909);
\coordinate (b27) at (-0.41215,0.22492);
\coordinate (b28) at (-0.35369,0.14610);
\coordinate (b29) at (-0.28097,0.08020);
\coordinate (b30) at (-0.19226,0.03824);
\coordinate (b31) at (-0.09707,0.01440);
\draw[outline] (b0) -- (b1) -- (b2) -- (b3) -- (b4) -- (b5) -- (b6) -- (b7) -- (b8) -- (b9) -- (b10) -- (b11) -- (b12) -- (b13) -- (b14) -- (b15) -- (b16) -- (b17) -- (b18) -- (b19) -- (b20) -- (b21) -- (b22) -- (b23) -- (b24) -- (b25) -- (b26) -- (b27) -- (b28) -- (b29) -- (b30) -- (b31) -- cycle;
\draw[pend] (b1) -- (b17);
\draw[pend] (b3) -- (b19);
\draw[pend] (b5) -- (b21);
\draw[pend] (b7) -- (b23);
\draw[pend] (b9) -- (b25);
\draw[pend] (b11) -- (b27);
\draw[pend] (b13) -- (b29);
\draw[pend] (b15) -- (b31);
\draw[path] (b0) -- (b15);
\draw[path] (b0) -- (b17);
\draw[path] (b2) -- (b17);
\draw[path] (b2) -- (b18);
\draw[path] (b3) -- (b18);
\draw[path] (b3) -- (b20);
\draw[path] (b4) -- (b20);
\draw[path] (b4) -- (b21);
\draw[path] (b6) -- (b21);
\draw[path] (b6) -- (b22);
\draw[path] (b7) -- (b22);
\draw[path] (b7) -- (b24);
\draw[path] (b8) -- (b25);
\draw[path] (b10) -- (b25);
\draw[path] (b10) -- (b26);
\draw[path] (b11) -- (b26);
\draw[path] (b11) -- (b28);
\draw[path] (b12) -- (b28);
\draw[path] (b12) -- (b29);
\draw[path] (b14) -- (b29);
\draw[path] (b14) -- (b30);
\draw[path] (b15) -- (b30);
\draw[axis] (b0) -- (b16);
\end{tikzpicture}

\small (b) \(A_{32}\)
\end{minipage}\hfill
\begin{minipage}[b]{0.32\textwidth}\centering
\begin{tikzpicture}[scale=3.8,
   outline/.style={dashed,line width=0.4pt},
   axis/.style={red!80!black,densely dotted,line width=0.6pt},
   path/.style={black,line width=0.3pt},
   pend/.style={blue!75!black,line width=0.6pt},
   lbl/.style={font=\tiny,inner sep=1pt}]
\coordinate (c0) at (0.00000,0.00000);
\coordinate (c1) at (0.04895,0.00361);
\coordinate (c2) at (0.09766,0.00962);
\coordinate (c3) at (0.14555,0.02037);
\coordinate (c4) at (0.19135,0.03804);
\coordinate (c5) at (0.23519,0.06011);
\coordinate (c6) at (0.27789,0.08430);
\coordinate (c7) at (0.31802,0.11256);
\coordinate (c8) at (0.35357,0.14640);
\coordinate (c9) at (0.38563,0.18357);
\coordinate (c10) at (0.41583,0.22226);
\coordinate (c11) at (0.44209,0.26373);
\coordinate (c12) at (0.46198,0.30860);
\coordinate (c13) at (0.47737,0.35521);
\coordinate (c14) at (0.49046,0.40251);
\coordinate (c15) at (0.49885,0.45087);
\coordinate (c16) at (0.49994,0.49994);
\coordinate (c17) at (0.49633,0.54889);
\coordinate (c18) at (0.49032,0.59760);
\coordinate (c19) at (0.47957,0.64549);
\coordinate (c20) at (0.46190,0.69129);
\coordinate (c21) at (0.43983,0.73513);
\coordinate (c22) at (0.41564,0.77784);
\coordinate (c23) at (0.38738,0.81796);
\coordinate (c24) at (0.35354,0.85351);
\coordinate (c25) at (0.31637,0.88557);
\coordinate (c26) at (0.27768,0.91577);
\coordinate (c27) at (0.23621,0.94203);
\coordinate (c28) at (0.19134,0.96192);
\coordinate (c29) at (0.14473,0.97731);
\coordinate (c30) at (0.09743,0.99040);
\coordinate (c31) at (0.04907,0.99880);
\coordinate (c32) at (0.00000,1.00000);
\coordinate (c33) at (-0.04907,0.99880);
\coordinate (c34) at (-0.09743,0.99040);
\coordinate (c35) at (-0.14473,0.97731);
\coordinate (c36) at (-0.19134,0.96192);
\coordinate (c37) at (-0.23621,0.94203);
\coordinate (c38) at (-0.27768,0.91577);
\coordinate (c39) at (-0.31637,0.88557);
\coordinate (c40) at (-0.35354,0.85351);
\coordinate (c41) at (-0.38738,0.81796);
\coordinate (c42) at (-0.41564,0.77784);
\coordinate (c43) at (-0.43983,0.73513);
\coordinate (c44) at (-0.46190,0.69129);
\coordinate (c45) at (-0.47957,0.64549);
\coordinate (c46) at (-0.49032,0.59760);
\coordinate (c47) at (-0.49633,0.54889);
\coordinate (c48) at (-0.49994,0.49994);
\coordinate (c49) at (-0.49885,0.45087);
\coordinate (c50) at (-0.49046,0.40251);
\coordinate (c51) at (-0.47737,0.35521);
\coordinate (c52) at (-0.46198,0.30860);
\coordinate (c53) at (-0.44209,0.26373);
\coordinate (c54) at (-0.41583,0.22226);
\coordinate (c55) at (-0.38563,0.18357);
\coordinate (c56) at (-0.35357,0.14640);
\coordinate (c57) at (-0.31802,0.11256);
\coordinate (c58) at (-0.27789,0.08430);
\coordinate (c59) at (-0.23519,0.06011);
\coordinate (c60) at (-0.19135,0.03804);
\coordinate (c61) at (-0.14555,0.02037);
\coordinate (c62) at (-0.09766,0.00962);
\coordinate (c63) at (-0.04895,0.00361);
\draw[outline] (c0) -- (c1) -- (c2) -- (c3) -- (c4) -- (c5) -- (c6) -- (c7) -- (c8) -- (c9) -- (c10) -- (c11) -- (c12) -- (c13) -- (c14) -- (c15) -- (c16) -- (c17) -- (c18) -- (c19) -- (c20) -- (c21) -- (c22) -- (c23) -- (c24) -- (c25) -- (c26) -- (c27) -- (c28) -- (c29) -- (c30) -- (c31) -- (c32) -- (c33) -- (c34) -- (c35) -- (c36) -- (c37) -- (c38) -- (c39) -- (c40) -- (c41) -- (c42) -- (c43) -- (c44) -- (c45) -- (c46) -- (c47) -- (c48) -- (c49) -- (c50) -- (c51) -- (c52) -- (c53) -- (c54) -- (c55) -- (c56) -- (c57) -- (c58) -- (c59) -- (c60) -- (c61) -- (c62) -- (c63) -- cycle;
\draw[pend] (c1) -- (c33);
\draw[pend] (c3) -- (c35);
\draw[pend] (c5) -- (c37);
\draw[pend] (c7) -- (c39);
\draw[pend] (c9) -- (c41);
\draw[pend] (c11) -- (c43);
\draw[pend] (c13) -- (c45);
\draw[pend] (c15) -- (c47);
\draw[pend] (c17) -- (c49);
\draw[pend] (c19) -- (c51);
\draw[pend] (c21) -- (c53);
\draw[pend] (c23) -- (c55);
\draw[pend] (c25) -- (c57);
\draw[pend] (c27) -- (c59);
\draw[pend] (c29) -- (c61);
\draw[pend] (c31) -- (c63);
\draw[path] (c0) -- (c31);
\draw[path] (c0) -- (c33);
\draw[path] (c2) -- (c33);
\draw[path] (c2) -- (c34);
\draw[path] (c3) -- (c34);
\draw[path] (c3) -- (c36);
\draw[path] (c4) -- (c36);
\draw[path] (c4) -- (c37);
\draw[path] (c6) -- (c37);
\draw[path] (c6) -- (c38);
\draw[path] (c7) -- (c38);
\draw[path] (c7) -- (c40);
\draw[path] (c8) -- (c40);
\draw[path] (c8) -- (c41);
\draw[path] (c10) -- (c41);
\draw[path] (c10) -- (c42);
\draw[path] (c11) -- (c42);
\draw[path] (c11) -- (c44);
\draw[path] (c12) -- (c44);
\draw[path] (c12) -- (c45);
\draw[path] (c14) -- (c45);
\draw[path] (c14) -- (c46);
\draw[path] (c15) -- (c46);
\draw[path] (c15) -- (c48);
\draw[path] (c16) -- (c49);
\draw[path] (c18) -- (c49);
\draw[path] (c18) -- (c50);
\draw[path] (c19) -- (c50);
\draw[path] (c19) -- (c52);
\draw[path] (c20) -- (c52);
\draw[path] (c20) -- (c53);
\draw[path] (c22) -- (c53);
\draw[path] (c22) -- (c54);
\draw[path] (c23) -- (c54);
\draw[path] (c23) -- (c56);
\draw[path] (c24) -- (c56);
\draw[path] (c24) -- (c57);
\draw[path] (c26) -- (c57);
\draw[path] (c26) -- (c58);
\draw[path] (c27) -- (c58);
\draw[path] (c27) -- (c60);
\draw[path] (c28) -- (c60);
\draw[path] (c28) -- (c61);
\draw[path] (c30) -- (c61);
\draw[path] (c30) -- (c62);
\draw[path] (c31) -- (c62);
\draw[axis] (c0) -- (c32);
\end{tikzpicture}

\small (c) \(A_{64}\)
\end{minipage}
\caption{The polygons \(A_n\),
numbered anticlockwise from \(p_0=(0,0)\) at the bottom (labels shown for
\(A_{16}\)), with the diameter graph: the
axis \(\{0,n/2\}\) (red, dotted), the pendant edges (blue), and the path from
\(0\) to \(n/4\) with its mirror image (thin). Dashed lines outline the polygons. Every other vertex pair is at distance at most
\(1-\delta(A_n)<1\), with \(\delta(A_n)\) the inactive
diameter margin of Definition~\ref{def:inactive-margin}
(Proposition~\ref{prop:a16-certified-checks}(4) and
Table~\ref{tab:family}).}
\label{fig:a16}
\end{figure}

\begin{figure}[t]
\centering
\begin{minipage}[b]{0.48\textwidth}\centering
\begin{tikzpicture}[scale=3.8,
   outline/.style={dashed,line width=0.4pt},
   axis/.style={red!80!black,densely dotted,line width=0.6pt},
   path/.style={black,line width=0.3pt},
   pend/.style={blue!75!black,line width=0.6pt}]
\coordinate (z0) at (0.00000,0.00000);
\coordinate (z1) at (0.08417,0.05045);
\coordinate (z2) at (0.16300,0.10891);
\coordinate (z3) at (0.23571,0.17481);
\coordinate (z4) at (0.30161,0.24753);
\coordinate (z5) at (0.36007,0.32635);
\coordinate (z6) at (0.41052,0.41052);
\coordinate (z7) at (0.45248,0.49923);
\coordinate (z8) at (0.48554,0.59163);
\coordinate (z9) at (0.49994,0.68871);
\coordinate (z10) at (0.49524,0.78673);
\coordinate (z11) at (0.47140,0.88192);
\coordinate (z12) at (0.38268,0.92388);
\coordinate (z13) at (0.29028,0.95694);
\coordinate (z14) at (0.19509,0.98079);
\coordinate (z15) at (0.09802,0.99518);
\coordinate (z16) at (-0.00000,1.00000);
\coordinate (z17) at (-0.09802,0.99518);
\coordinate (z18) at (-0.19509,0.98079);
\coordinate (z19) at (-0.29028,0.95694);
\coordinate (z20) at (-0.38268,0.92388);
\coordinate (z21) at (-0.47140,0.88192);
\coordinate (z22) at (-0.49524,0.78673);
\coordinate (z23) at (-0.49994,0.68871);
\coordinate (z24) at (-0.48554,0.59163);
\coordinate (z25) at (-0.45248,0.49923);
\coordinate (z26) at (-0.41052,0.41052);
\coordinate (z27) at (-0.36007,0.32635);
\coordinate (z28) at (-0.30161,0.24753);
\coordinate (z29) at (-0.23571,0.17481);
\coordinate (z30) at (-0.16300,0.10891);
\coordinate (z31) at (-0.08417,0.05045);
\draw[outline] (z0) -- (z1) -- (z2) -- (z3) -- (z4) -- (z5) -- (z6) -- (z7) -- (z8) -- (z9) -- (z10) -- (z11) -- (z12) -- (z13) -- (z14) -- (z15) -- (z16) -- (z17) -- (z18) -- (z19) -- (z20) -- (z21) -- (z22) -- (z23) -- (z24) -- (z25) -- (z26) -- (z27) -- (z28) -- (z29) -- (z30) -- (z31) -- cycle;
\draw[pend] (z0) -- (z12);
\draw[pend] (z0) -- (z13);
\draw[pend] (z0) -- (z14);
\draw[pend] (z0) -- (z15);
\draw[pend] (z11) -- (z25);
\draw[pend] (z11) -- (z26);
\draw[pend] (z11) -- (z27);
\draw[pend] (z11) -- (z28);
\draw[pend] (z11) -- (z29);
\draw[pend] (z11) -- (z30);
\draw[pend] (z11) -- (z31);
\draw[pend] (z0) -- (z20);
\draw[pend] (z0) -- (z19);
\draw[pend] (z0) -- (z18);
\draw[pend] (z0) -- (z17);
\draw[pend] (z21) -- (z7);
\draw[pend] (z21) -- (z6);
\draw[pend] (z21) -- (z5);
\draw[pend] (z21) -- (z4);
\draw[pend] (z21) -- (z3);
\draw[pend] (z21) -- (z2);
\draw[pend] (z21) -- (z1);
\draw[path] (z0) -- (z11);
\draw[path] (z11) -- (z24);
\draw[path] (z24) -- (z10);
\draw[path] (z10) -- (z23);
\draw[path] (z0) -- (z21);
\draw[path] (z21) -- (z8);
\draw[path] (z8) -- (z22);
\draw[path] (z22) -- (z9);
\draw[axis] (z0) -- (z16);
\end{tikzpicture}

\small $Z_{32}$\end{minipage}\hfill%
\begin{minipage}[b]{0.48\textwidth}\centering
\begin{tikzpicture}[scale=3.8,
   outline/.style={dashed,line width=0.4pt},
   axis/.style={red!80!black,densely dotted,line width=0.6pt},
   path/.style={black,line width=0.3pt},
   pend/.style={blue!75!black,line width=0.6pt}]
\coordinate (w0) at (0.00000,0.00000);
\coordinate (w1) at (0.04895,0.00361);
\coordinate (w2) at (0.09731,0.01200);
\coordinate (w3) at (0.14520,0.02276);
\coordinate (w4) at (0.19180,0.03815);
\coordinate (w5) at (0.23668,0.05804);
\coordinate (w6) at (0.28052,0.08011);
\coordinate (w7) at (0.32198,0.10637);
\coordinate (w8) at (0.36068,0.13657);
\coordinate (w9) at (0.39623,0.17041);
\coordinate (w10) at (0.42829,0.20758);
\coordinate (w11) at (0.45655,0.24771);
\coordinate (w12) at (0.47862,0.29155);
\coordinate (w13) at (0.49401,0.33815);
\coordinate (w14) at (0.50000,0.38687);
\coordinate (w15) at (0.49880,0.43594);
\coordinate (w16) at (0.49516,0.48488);
\coordinate (w17) at (0.48677,0.53324);
\coordinate (w18) at (0.47602,0.58113);
\coordinate (w19) at (0.46293,0.62844);
\coordinate (w20) at (0.44526,0.67423);
\coordinate (w21) at (0.42537,0.71910);
\coordinate (w22) at (0.40118,0.76181);
\coordinate (w23) at (0.37492,0.80328);
\coordinate (w24) at (0.34472,0.84197);
\coordinate (w25) at (0.31088,0.87752);
\coordinate (w26) at (0.27371,0.90958);
\coordinate (w27) at (0.23358,0.93784);
\coordinate (w28) at (0.19088,0.96203);
\coordinate (w29) at (0.14509,0.97970);
\coordinate (w30) at (0.09778,0.99279);
\coordinate (w31) at (0.04907,0.99880);
\coordinate (w32) at (-0.00000,1.00000);
\coordinate (w33) at (-0.04907,0.99880);
\coordinate (w34) at (-0.09778,0.99279);
\coordinate (w35) at (-0.14509,0.97970);
\coordinate (w36) at (-0.19088,0.96203);
\coordinate (w37) at (-0.23358,0.93784);
\coordinate (w38) at (-0.27371,0.90958);
\coordinate (w39) at (-0.31088,0.87752);
\coordinate (w40) at (-0.34472,0.84197);
\coordinate (w41) at (-0.37492,0.80328);
\coordinate (w42) at (-0.40118,0.76181);
\coordinate (w43) at (-0.42537,0.71910);
\coordinate (w44) at (-0.44526,0.67423);
\coordinate (w45) at (-0.46293,0.62844);
\coordinate (w46) at (-0.47602,0.58113);
\coordinate (w47) at (-0.48677,0.53324);
\coordinate (w48) at (-0.49516,0.48488);
\coordinate (w49) at (-0.49880,0.43594);
\coordinate (w50) at (-0.50000,0.38687);
\coordinate (w51) at (-0.49401,0.33815);
\coordinate (w52) at (-0.47862,0.29155);
\coordinate (w53) at (-0.45655,0.24771);
\coordinate (w54) at (-0.42829,0.20758);
\coordinate (w55) at (-0.39623,0.17041);
\coordinate (w56) at (-0.36068,0.13657);
\coordinate (w57) at (-0.32198,0.10637);
\coordinate (w58) at (-0.28052,0.08011);
\coordinate (w59) at (-0.23668,0.05804);
\coordinate (w60) at (-0.19180,0.03815);
\coordinate (w61) at (-0.14520,0.02276);
\coordinate (w62) at (-0.09731,0.01200);
\coordinate (w63) at (-0.04895,0.00361);
\draw[outline] (w0) -- (w1) -- (w2) -- (w3) -- (w4) -- (w5) -- (w6) -- (w7) -- (w8) -- (w9) -- (w10) -- (w11) -- (w12) -- (w13) -- (w14) -- (w15) -- (w16) -- (w17) -- (w18) -- (w19) -- (w20) -- (w21) -- (w22) -- (w23) -- (w24) -- (w25) -- (w26) -- (w27) -- (w28) -- (w29) -- (w30) -- (w31) -- (w32) -- (w33) -- (w34) -- (w35) -- (w36) -- (w37) -- (w38) -- (w39) -- (w40) -- (w41) -- (w42) -- (w43) -- (w44) -- (w45) -- (w46) -- (w47) -- (w48) -- (w49) -- (w50) -- (w51) -- (w52) -- (w53) -- (w54) -- (w55) -- (w56) -- (w57) -- (w58) -- (w59) -- (w60) -- (w61) -- (w62) -- (w63) -- cycle;
\draw[pend] (w30) -- (w62);
\draw[pend] (w28) -- (w59);
\draw[pend] (w53) -- (w22);
\draw[pend] (w52) -- (w20);
\draw[pend] (w51) -- (w18);
\draw[pend] (w51) -- (w17);
\draw[pend] (w34) -- (w2);
\draw[pend] (w36) -- (w5);
\draw[pend] (w11) -- (w42);
\draw[pend] (w12) -- (w44);
\draw[pend] (w13) -- (w46);
\draw[pend] (w13) -- (w47);
\draw[path] (w0) -- (w31);
\draw[path] (w31) -- (w63);
\draw[path] (w63) -- (w30);
\draw[path] (w30) -- (w61);
\draw[path] (w61) -- (w29);
\draw[path] (w29) -- (w60);
\draw[path] (w60) -- (w28);
\draw[path] (w28) -- (w58);
\draw[path] (w58) -- (w27);
\draw[path] (w27) -- (w57);
\draw[path] (w57) -- (w26);
\draw[path] (w26) -- (w56);
\draw[path] (w56) -- (w25);
\draw[path] (w25) -- (w55);
\draw[path] (w55) -- (w24);
\draw[path] (w24) -- (w54);
\draw[path] (w54) -- (w23);
\draw[path] (w23) -- (w53);
\draw[path] (w53) -- (w21);
\draw[path] (w21) -- (w52);
\draw[path] (w52) -- (w19);
\draw[path] (w19) -- (w51);
\draw[path] (w51) -- (w16);
\draw[path] (w15) -- (w50);
\draw[path] (w50) -- (w14);
\draw[path] (w14) -- (w49);
\draw[path] (w0) -- (w33);
\draw[path] (w33) -- (w1);
\draw[path] (w1) -- (w34);
\draw[path] (w34) -- (w3);
\draw[path] (w3) -- (w35);
\draw[path] (w35) -- (w4);
\draw[path] (w4) -- (w36);
\draw[path] (w36) -- (w6);
\draw[path] (w6) -- (w37);
\draw[path] (w37) -- (w7);
\draw[path] (w7) -- (w38);
\draw[path] (w38) -- (w8);
\draw[path] (w8) -- (w39);
\draw[path] (w39) -- (w9);
\draw[path] (w9) -- (w40);
\draw[path] (w40) -- (w10);
\draw[path] (w10) -- (w41);
\draw[path] (w41) -- (w11);
\draw[path] (w11) -- (w43);
\draw[path] (w43) -- (w12);
\draw[path] (w12) -- (w45);
\draw[path] (w45) -- (w13);
\draw[path] (w13) -- (w48);
\draw[axis] (w0) -- (w32);
\end{tikzpicture}

\small $Z_{64}$\end{minipage}
\caption{The reconstructed polygons $Z_{32}$ and $Z_{64}$, with their diameter graphs. The axis is red and dotted, the paths are thin black, the pendant edges are blue, and the polygon boundaries are dashed.}
\label{fig:zn}
\end{figure}
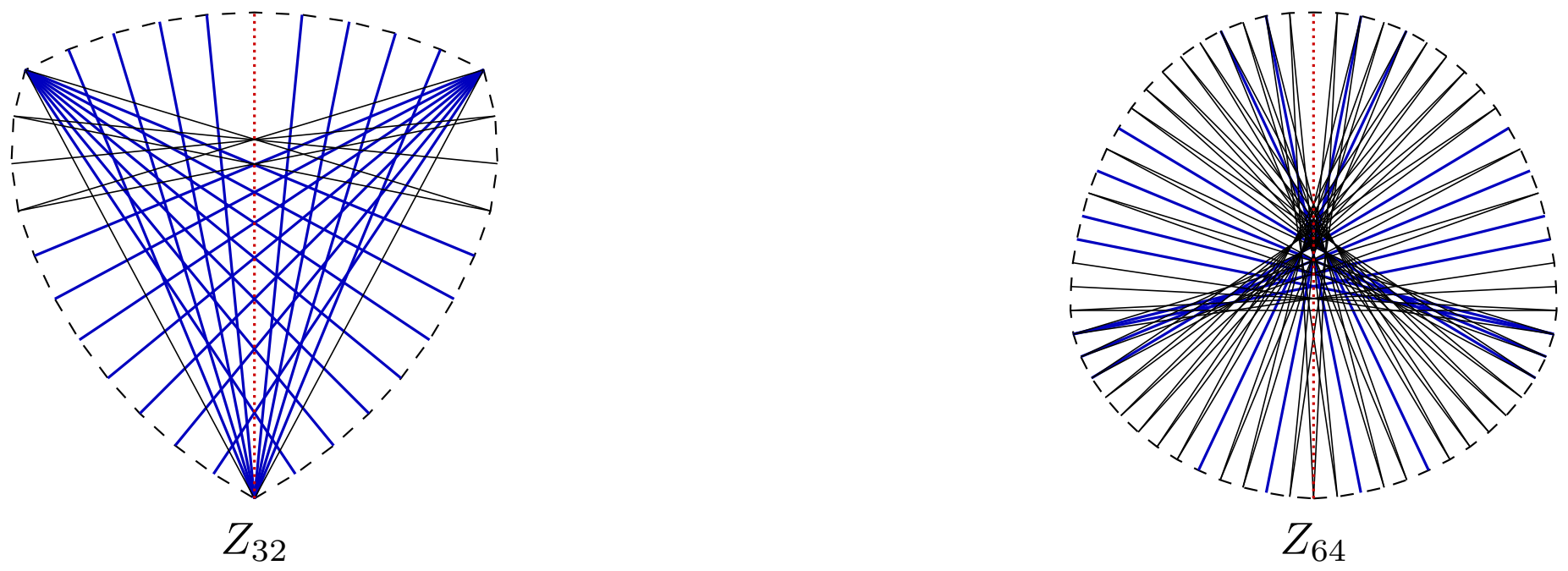

We illustrate the criterion with the square which satisfies the
positive multiplier and rigidity conditions, and a rhombus admitting
an explicit perimeter increasing deformation. Our main
applications concern the polygons $A_{16}$, $A_{32}$, and $A_{64}$
obtained from the construction of Bingane and
Audet~\cite{binganeaudet2022}, together with reconstructions of
$Z_{32}$ and $Z_{64}$ from their Propositions~1 and~2.
Figure~\ref{fig:a16} shows the diameter graphs of $A_{16}$,
$A_{32}$ and $A_{64}$, while Figure~\ref{fig:zn} shows those of
$Z_{32}$ and $Z_{64}$. Using exact algebra and validated interval
arithmetic, we establish the active constraints, strictly positive
multipliers, and infinitesimal rigidity for all five configurations.
Our criterion then shows that each is a strict local maximiser. Moreover, since
$L(A_n)<L(Z_n)$ for $n\in\{32,64\}$,
we conclude that $A_{32}$ and $A_{64}$ are not global maximisers. Therefore,
strict local maximisers with different perimeters coexist in the space of
convex equilateral small polygons with the same number of sides. An accompanying Lean~4 formalisation independently verifies
the criterion and the local maximality results.

Although our applications concern perimeter, the admissible
cone and rigidity conditions depend only on the polygonal
constraints. The same calculus therefore applies to other smooth
Euclidean invariant functionals whenever the corresponding
positive multiplier balance holds.

The paper is organised as follows. Section~2 sets up the
configuration space, active constraints and associated linearised
spaces. Section~3 proves the first order criterion and its
instability counterpart. Section~4 gives the square and rhombus
examples. Section~5 describes the $A_n$ construction and details
$A_{16}$. Section~6 establishes strict local maximality for
$A_{16}$, $A_{32}$, $A_{64}$, $Z_{32}$ and $Z_{64}$, and gives
an explicit radius for $A_{16}$. The ESP1 software package~\cite{chiqpkg} provides the computational
data, verification code, and Lean~4 formalisation.

\section{Configuration Space and Realised Diameter Geometry}
\label{sec:first-order-setup}
Throughout, $\langle\cdot,\cdot\rangle$ denotes the Euclidean
inner product on $\R^2$, extended to $(\R^2)^n$ by
$\langle v,w\rangle=\sum_i\langle v_i,w_i\rangle$ and to
$(\R^2)^n\times\R$ by
$\langle(v,\eta),(w,\zeta)\rangle
=\langle v,w\rangle+\eta\zeta$.
We use $|\cdot|$ for the Euclidean norm on $\R^2$
and $\|\cdot\|$ for the induced norms on the product spaces.
Orthogonal complements are taken in the relevant ambient
inner product space. Vertex indices are cyclic modulo $n$.

\begin{definition}
\label{def:small-polygon}
A small $n$-gon is a labelled cyclic tuple
$p=(p_0,\ldots,p_{n-1})\in(\R^2)^n$
satisfying
\[
|p_i-p_j|\leq 1 \qquad (i\neq j).
\]
We call $p$ strictly convex if its vertices are pairwise
distinct, no three are collinear, every vertex is an
extreme point of their convex hull, and the labelling
follows the positive cyclic order.
\end{definition}

Note that no convexity or geometric cyclic ordering of the vertices
is required in this definition.

\begin{definition}
An equilateral small $n$-gon is a pair $(p,\ell)$, where
$p$ is a small $n$-gon in the sense of
Definition~\ref{def:small-polygon} and $\ell>0$ satisfies
\[
|p_{i+1}-p_i|=\ell
\qquad (i=0,\ldots,n-1).
\]
Its perimeter is $L(p,\ell)=n\ell$, and the pair
is strictly convex when $p$ is strictly convex.

We call $(p,\ell)$ feasible if $p$ belongs to the
closure in $(\R^2)^n$ of the strictly convex tuples
labelled in positive cyclic order. Feasible pairs
are also called convex equilateral small $n$-gons.
\end{definition}

\begin{definition}
The set of active diameter edges of \(p\) is
\(E(p) = \bigl\{\{i,j\}:|p_i-p_j|=1\bigr\}\). The graph on the vertices of
\(p\) with edge set \(E(p)\) is the \emph{diameter graph} of \(p\).
\end{definition}

\begin{definition}
\label{def:inactive-margin}
For a small $n$-gon $p$, with $n\geq4$, define
the minimum inactive diameter slack by
\[
\delta(p)
=
\min_{\{i,j\}\notin E(p)}
\bigl(1-|p_i-p_j|\bigr).
\]
\end{definition}

Fix an equilateral small \(n\)-gon \((p,\ell)\). For each side let
\(s_i=(p_{i+1}-p_i)/\ell\) be its unit vector, and for each active diameter
edge let \(u_{ij}=p_j-p_i\), so that \(|u_{ij}|=1\).

The side constraint functions are \(e_i(q,m)=|q_{i+1}-q_i|-m\) for \(i=0,\ldots,n-1\), and the
diameter constraint functions are \(d_{ij}(q)=|q_i-q_j|-1\) for every pair \(\{i,j\}\); a
competitor \((q,m)\) is small exactly when \(d_{ij}(q)\leq0\) for all pairs. Since
\(p_{i+1}\neq p_i\) and \(p_i\neq p_j\) for \(\{i,j\}\in E(p)\), \(e_i\) is differentiable at
\((p,\ell)\) and \(d_{ij}\) at \(p\), with
\begin{equation}
\label{eq:constraint-differentials}
De_i(p,\ell)[v,\eta]
=
\langle s_i,v_{i+1}-v_i\rangle-\eta,
\qquad
Dd_{ij}(p)[v]
=
\langle u_{ij},v_j-v_i\rangle .
\end{equation}
Only the constraints of active pairs, for which \(d_{ij}(p)=0\), matter to first
order: an active diameter may break inward, or stay active, precisely when
\(Dd_{ij}(p)[v]\leq0\). Swapping \(i\) and \(j\) reverses the signs of both
\(u_{ij}\) and \(v_j-v_i\), so the product does not depend on how an edge is oriented.

\begin{definition}
\label{def:cone}
The linearised feasible cone at \((p,\ell)\) is
\[
\Kcal(p,\ell)
=
\left\{
(v,\eta)\in(\R^2)^n\times\R:
\begin{array}{l}
\langle s_i,v_{i+1}-v_i\rangle=\eta,
\quad i=0,\ldots,n-1,\\[2mm]
\langle u_{ij},v_j-v_i\rangle\leq0,
\quad \{i,j\}\in E(p)
\end{array}
\right\}.
\]
Membership does not mean that a direction is realised by a feasible curve; see
Proposition~\ref{prop:instability}.
\end{definition}
\begin{definition}
\label{def:equality}
The active equality space is
\[
\Ical(p,\ell)
=
\left\{
(v,\eta)\in(\R^2)^n\times\R:
\begin{array}{l}
\langle s_i,v_{i+1}-v_i\rangle=\eta,
\quad i=0,\ldots,n-1,\\[2mm]
\langle u_{ij},v_j-v_i\rangle=0,
\quad \{i,j\}\in E(p)
\end{array}
\right\}.
\]
Every direction in \(\Ical(p,\ell)\) satisfies the active diameter
conditions of Definition~\ref{def:cone} with equality, so
\[
\Ical(p,\ell)\subseteq\Kcal(p,\ell).
\]
\end{definition}

Let
\(J(x,y)=(-y,x)\).
The space of infinitesimal Euclidean motions at \(p\) is \[
R(p)=\left\{(a+\omega Jp_0,\ldots,a+\omega Jp_{n-1}) : a\in\R^2,\ \omega\in\R\right\}.
\] Inside \((\R^2)^n\times\R\), set
\(R_{\mathrm{eq}}(p) = \{(r,0):r\in R(p)\}\).
Then \(R_{\mathrm{eq}}(p)\subseteq\Ical(p,\ell)\), and
\(\dim R_{\mathrm{eq}}(p)=3\): if \(a+\omega Jp_i=0\) for every \(i\), subtracting at the
distinct vertices \(p_0\neq p_1\) (as \(\ell>0\)) gives \(\omega J(p_0-p_1)=0\), so \(\omega=0\) and then \(a=0\).

\begin{definition}
\label{def:jacobian}
Write \(B=D(e,d)\) for the Jacobian at \((p,\ell)\) of the constraint map \[
(q,m)\longmapsto\left(e_0(q,m),\ldots,e_{n-1}(q,m),(d_{ij}(q))_{\{i,j\}\in E(p)}\right).
\] By \eqref{eq:constraint-differentials} its kernel is the active equality space
\(\Ical(p,\ell)\) of Definition~\ref{def:equality}. Dependent rows are allowed. As the ambient space has dimension \(2n+1\) and
the kernel contains \(R_{\mathrm{eq}}(p)\), \(\rank B\leq2n-2\), with equality if and only if
\(\Ical(p,\ell)=R_{\mathrm{eq}}(p)\).
\end{definition}

\section{First Order Criteria for Local Maximality}
\label{sec:first-order}
Since \(L(p,\ell)=n\ell\) is linear in \(\ell\) and independent of \(p\),
\begin{equation}
\label{eq:first-variation}
DL(p,\ell)[v,\eta]=n\eta .
\end{equation}
Given side multipliers \(\alpha_i\in\R\) and active diameter multipliers
\(\lambda_{ij}\in\R\), we say that the \emph{multiplier identity} holds if
\begin{equation}
\label{eq:multiplier-identity}
DL(p,\ell)[v,\eta]
=
\sum_{i=0}^{n-1}
\alpha_i\,De_i(p,\ell)[v,\eta]
+
\sum_{\{i,j\}\in E(p)}
\lambda_{ij}\,Dd_{ij}(p)[v]
\end{equation}
for every \((v,\eta)\in(\R^2)^n\times\R\).
\begin{proposition}
\label{prop:first-order-no-ascent}
Let \((p,\ell)\) be an equilateral small \(n\)-gon. Suppose
\eqref{eq:multiplier-identity} holds for some side multipliers
\(\alpha_i\in\R\) and some active diameter multipliers
\(\lambda_{ij}\geq0 \quad (\{i,j\}\in E(p))\).
Then
\(DL(p,\ell)[v,\eta]\leq0\)
for every
\((v,\eta)\in\Kcal(p,\ell)\).
\end{proposition}

\begin{proof}
For a direction in \(\Kcal(p,\ell)\), the two formulas
\eqref{eq:constraint-differentials} show that all side differentials vanish and
that every active diameter differential is nonpositive. Hence
\[
DL(p,\ell)[v,\eta] = \sum_{\{i,j\}\in E(p)} \lambda_{ij} \langle u_{ij},v_j-v_i\rangle \leq0.
\qedhere
\]
\end{proof}

\begin{remark}
\label{rem:multiplier-balance}
Comparison of the \(\eta\) coefficients in \eqref{eq:multiplier-identity}, using
\eqref{eq:first-variation}, gives
\(\sum_{i=0}^{n-1}\alpha_i=-n\), providing a consistency check
for computed multipliers.
\end{remark}

\begin{lemma}
\label{lem:positive-multipliers}
Let \((p,\ell)\) be an equilateral small \(n\)-gon and assume
\eqref{eq:multiplier-identity} holds with
\(\lambda_{ij}>0 \quad (\{i,j\}\in E(p))\).
If
\((v,\eta)\in\Kcal(p,\ell)\)
and
\(DL(p,\ell)[v,\eta]=0\),
then
\((v,\eta)\in\Ical(p,\ell)\).
\end{lemma}

\begin{proof}
As in the proof of Proposition~\ref{prop:first-order-no-ascent},
\(DL(p,\ell)[v,\eta]=\sum_{\{i,j\}\in E(p)}\lambda_{ij}\langle u_{ij},v_j-v_i\rangle\). Every summand is nonpositive. Since every multiplier is strictly positive, a
zero sum forces
\(\langle u_{ij},v_j-v_i\rangle=0\)
for every active diameter edge.
\end{proof}

\subsection{A local Euclidean slice}
\begin{lemma}
\label{lem:euclidean-slice}
Let
\(p=(p_0,\ldots,p_{n-1})\in(\R^2)^n\)
have at least two distinct vertices. There are a neighbourhood \(\mathcal U\) of
\(p\) and a neighbourhood \(\mathcal V\) of the identity in the group \(\mathrm{SE}(2)\) of
orientation preserving Euclidean motions, acting on tuples vertexwise, such that, for every \(q\in\mathcal U\), there is a unique
\(\gamma(q)\in\mathcal V\) satisfying
\(\gamma(q)q-p\in R(p)^\perp\).
The map
\(q\longmapsto \gamma(q)\)
is smooth and
\(\gamma(p)=\operatorname{id}\).
\end{lemma}

\begin{proof}
Let \(\tau_1=((1,0),\ldots,(1,0))\), \(\tau_2=((0,1),\ldots,(0,1))\) and
\(\tau_3=(Jp_0,\ldots,Jp_{n-1})\).
These vectors form a basis of \(R(p)\): if \(a+\omega Jp_i=0\) for every \(i\), then \(\omega J(p_i-p_j)=0\) at two distinct vertices, so \(\omega=0\) and \(a=0\).
Write an orientation preserving Euclidean motion near the identity as
\(g_{a,\theta}(x)=e^{\theta J}x+a\), where \(e^{\theta J}\) is the rotation by \(\theta\).
Define
\[
\Phi(a,\theta,q)
=
\left(
\langle g_{a,\theta}q-p,\tau_1\rangle,
\langle g_{a,\theta}q-p,\tau_2\rangle,
\langle g_{a,\theta}q-p,\tau_3\rangle
\right).
\]

Then
\(\Phi(0,0,p)=0\).
The derivative with respect to \((a,\theta)\) at \((0,0,p)\) is the Gram
matrix of \(\tau_1,\tau_2,\tau_3\), and is therefore invertible. The implicit function
theorem gives unique smooth functions \(a(q)\) and \(\theta(q)\) satisfying
\(\Phi(a(q),\theta(q),q)=0\).
Setting
\(\gamma(q)=g_{a(q),\theta(q)}\)
gives the conclusion.
\end{proof}
\begin{figure}[ht]
\centering
\begin{tikzpicture}[scale=1.05,
   vtx/.style={circle,fill=black,inner sep=0pt,minimum size=2.6pt},
   gvtx/.style={circle,draw=black,fill=white,line width=0.4pt,inner sep=0pt,
                minimum size=2.8pt},
   lbl/.style={font=\small,inner sep=2pt},
   arr/.style={-latex,line width=0.55pt}]
\draw[line width=0.9pt] plot[smooth,tension=0.85]
  coordinates {(-2.35,-0.80) (-1.20,-0.24) (0,0) (1.20,-0.24) (2.35,-0.80)};
\node[lbl,anchor=west] at (2.42,-0.80) {orbit of \(p\)};
\draw[dashed,line width=0.5pt] (0,-0.55) -- (0,2.30);
\node[lbl,anchor=south] at (0,2.32) {\(p+R(p)^\perp\)};
\node[vtx] at (0,0) {}; \node[lbl,anchor=north east] at (-0.02,-0.04) {\(p\)};
\node[vtx] at (1.30,1.78) {}; \node[lbl,anchor=west] at (1.34,1.78) {\(q^1\)};
\node[vtx] at (0.82,1.12) {}; \node[lbl,anchor=west] at (0.86,1.12) {\(q^2\)};
\node[vtx] at (0.50,0.66) {}; \node[lbl,anchor=west] at (0.54,0.66) {\(q^3\)};
\node[gvtx] at (0,1.78) {};
\node[gvtx] at (0,1.12) {};
\node[gvtx] at (0,0.66) {};
\draw[arr,gray] (1.24,1.78) -- (0.10,1.78);
\draw[arr,gray] (0.76,1.12) -- (0.10,1.12);
\draw[arr,gray] (0.44,0.66) -- (0.10,0.66);
\draw[arr,line width=1.1pt] (0,0) -- (0,1.30);
\node[lbl,anchor=east] at (-0.06,0.95) {\(v\)};
\end{tikzpicture}
\caption{The proof of Theorem~\ref{thm:first-order-topological}, in configuration
space. Competitors
\(q^k\to p\) are moved by the slice of Lemma~\ref{lem:euclidean-slice} onto
\(p+R(p)^\perp\) (hollow). Their rescaled displacements accumulate on a unit
vector \((v,\eta)\in\Kcal(p,\ell)\) with \(\eta\geq0\) and
\(v\in R(p)^\perp\), which hypotheses (i) and (ii) rule out.}
\label{fig:slice}
\end{figure}
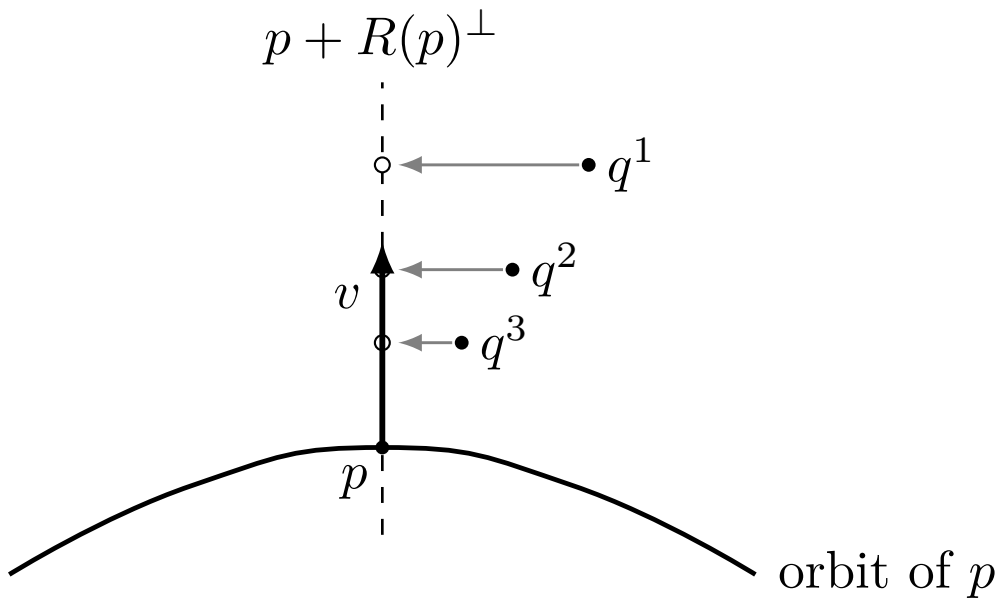

\begin{theorem}
\label{thm:first-order-topological}
Let \((p,\ell)\) be an equilateral small \(n\)-gon. Suppose:
\begin{enumerate}[label=\textup{(\roman*)}]
\item the multiplier identity \eqref{eq:multiplier-identity} holds with side
multipliers \(\alpha_i\in\R\) and strictly positive active diameter
multipliers;
\item \(\Ical(p,\ell)=R_{\mathrm{eq}}(p)\).
\end{enumerate}
Then \((p,\ell)\) is a strict local maximiser of perimeter modulo orientation preserving Euclidean motions among all equilateral small \(n\)-gons.
More precisely, there is a neighbourhood \(\mathcal W\) of \((p,\ell)\) such that every equilateral small \(n\)-gon \((q,m)\in\mathcal W\) satisfies $m\leq\ell$, with equality if and only if \(q\) is the image of \(p\) under an orientation-preserving Euclidean motion.
\end{theorem}

\begin{proof}
Figure~\ref{fig:slice} illustrates the argument. Suppose the conclusion fails. Then there are equilateral small
\(n\)-gons \((q^k,m_k)\to(p,\ell)\) with \(m_k\geq\ell\), none of which is
an image of \(p\). For large \(k\), \(q^k\) lies in the neighbourhood \(\mathcal U\) of
Lemma~\ref{lem:euclidean-slice}, and \(\gamma(q^k)\to\operatorname{id}\) since \(\gamma\)
is smooth with \(\gamma(p)=\operatorname{id}\). Replacing \(q^k\) by
\(\gamma(q^k)q^k\), which preserves the side and diameter constraints, convergence and the property of
not being an image of \(p\), we
may assume that
\[
q^k-p\in R(p)^\perp\qquad\text{for all }k.
\]

Set \(\varepsilon_k=\|(q^k-p,m_k-\ell)\|\), which is positive since \(q^k\neq p\). A subsequence of
\((q^k-p,m_k-\ell)/\varepsilon_k\) converges to a unit vector \((v,\eta)\). Dividing
\(|q^k_{i+1}-q^k_i|-|p_{i+1}-p_i|=m_k-\ell\) and, for active pairs,
\(|q^k_i-q^k_j|-|p_i-p_j|\leq0\) by \(\varepsilon_k\), and letting \(k\to\infty\), gives (the norm is differentiable at the nonzero vectors \(p_{i+1}-p_i\) and \(p_j-p_i\))
\[
\langle s_i,v_{i+1}-v_i\rangle=\eta,
\qquad
\langle u_{ij},v_j-v_i\rangle\leq0 .
\]
Thus \((v,\eta)\in\Kcal(p,\ell)\), and \(\eta\geq0\) because
\(m_k\geq\ell\).

By Proposition~\ref{prop:first-order-no-ascent} and \eqref{eq:first-variation},
\(n\eta=DL(p,\ell)[v,\eta]\leq0\), so \(\eta=0\).
Lemma~\ref{lem:positive-multipliers} and hypothesis~(ii) give
\((v,0)\in\Ical(p,\ell)=R_{\mathrm{eq}}(p)\), and \(v\in R(p)^\perp\) as a limit of
vectors in the closed subspace \(R(p)^\perp\); hence
\(v\in R(p)\cap R(p)^\perp=\{0\}\). This contradicts \(\|(v,\eta)\|=1\).
\end{proof}

\begin{remark}
Both hypotheses concern the single configuration \((p,\ell)\).
Hypothesis~(ii) is equivalent to \(\rank B=2n-2\)
(Definition~\ref{def:jacobian}) and beyond that nothing is
assumed about the constraints, which need not be independent,
nor about the smoothness of the constraint set or second
derivatives.
The proof uses only that \(p\) has two distinct vertices and
\(\ell>0\) and, of the competitors, only the side and active
diameter constraints.
The neighbourhood \(\mathcal W\) is not quantified; for
\(A_{16}\), Theorem~\ref{thm:a16-radius} gives an explicit radius.
\end{remark}

We conclude with a simple instability criterion.

\begin{proposition}
\label{prop:instability}
Let \((p,\ell)\) be a strictly convex equilateral small \(n\)-gon. Suppose that
there is a \(C^1\) curve \(\beta\) of feasible pairs with \(\beta(0)=(p,\ell)\)
whose initial velocity \((v,\eta)\) satisfies \(\eta>0\). Then \((p,\ell)\) is not a
local maximiser of perimeter among convex equilateral small \(n\)-gons.
\end{proposition}
\begin{proof}
Along \(\beta\), by \eqref{eq:first-variation}, \[
\frac{d}{dt}\bigg|_{t=0^+}L(\beta(t))=DL(p,\ell)[v,\eta]=n\eta>0.
\] Hence the perimeter increases for all sufficiently small positive \(t\).
\end{proof}

\section{Examples}
\label{sec:examples}
\subsection{The square}
 Let \(a=1/(2\sqrt2)\), and let \(p_0=(a,a)\), \(p_1=(-a,a)\), \(p_2=(-a,-a)\),
\(p_3=(a,-a)\) and \(\ell=1/\sqrt2\). The active diameter edges are
\(E(p) = \{\{0,2\},\{1,3\}\}\). The side unit vectors are \(s_0=(-1,0)\), \(s_1=(0,-1)\), \(s_2=(1,0)\), \(s_3=(0,1)\), while
\[
u_{02}
=
\left(
-\frac1{\sqrt2},
-\frac1{\sqrt2}
\right),
\qquad
u_{13}
=
\left(
\frac1{\sqrt2},
-\frac1{\sqrt2}
\right).
\]

\begin{lemma}
At the square,
\(\rank B=6\).
Moreover, \eqref{eq:multiplier-identity} holds with
\(\alpha_0=\alpha_1=\alpha_2=\alpha_3=-1\)
and \(\lambda_{02}=\lambda_{13}=\sqrt2\).
\end{lemma}
\begin{proof}
Suppose \(\sum_{i=0}^3\beta_iDe_i+\gamma_{02}Dd_{02}+\gamma_{13}Dd_{13}=0\). The \(\eta\) coefficient gives
\(\beta_0+\beta_1+\beta_2+\beta_3=0\).
The vertex coefficients give
\[
\begin{aligned}
-\beta_0s_0+\beta_3s_3-\gamma_{02}u_{02}&=0, &\quad \beta_0s_0-\beta_1s_1-\gamma_{13}u_{13}&=0,\\
\beta_1s_1-\beta_2s_2+\gamma_{02}u_{02}&=0, &\quad \beta_2s_2-\beta_3s_3+\gamma_{13}u_{13}&=0.
\end{aligned}
\]
Comparing components gives \(\beta_0=\beta_1=\beta_2=\beta_3=-\gamma_{02}/\sqrt2\) and
\(\gamma_{02}=\gamma_{13}\).
The \(\eta\) equation then forces
\(\gamma_{02}=0\).
Thus all six rows are independent and
\(\rank B=6\).

For \eqref{eq:multiplier-identity}, direct substitution of
\(\alpha_i=-1\)
and \(\lambda_{02}=\lambda_{13}=\sqrt2\) gives the required vertex balances, while \(-\sum_{i=0}^3\alpha_i=4\) matches the \(\eta\) coefficient of \(DL[v,\eta]=4\eta\).
\end{proof}
\begin{proposition}
The square satisfies \(\Ical(p,\ell)=R_{\mathrm{eq}}(p)\) and is a strict local maximiser of perimeter modulo orientation preserving
Euclidean motions among convex equilateral small quadrilaterals.
\end{proposition}
\begin{proof}
Since \(\rank B=6=2n-2\), the two spaces are equal by
Definition~\ref{def:jacobian}. The active diameter multipliers are strictly
positive, so Theorem~\ref{thm:first-order-topological} applies.
\end{proof}
\subsection{A nonoptimal rhombus}
Let
\[
p_0=(0,0),
\quad
p_1=\left(\frac12,-\frac25\right),
\quad
p_2=(1,0),
\quad
p_3=\left(\frac12,\frac25\right),
\]
and \(\ell=\sqrt{41}/10\). Then \(E(p)=\{\{0,2\}\}\), the other diagonal having length \(4/5\).
 \begin{proposition}
\label{prop:rhombus}
The rhombus is not a local maximiser of perimeter among convex equilateral small quadrilaterals.
\end{proposition}
\begin{proof}
Let \(v_0=v_2=(0,0)\), \(v_1=(0,-1)\), \(v_3=(0,1)\), and set \(\eta=2/(5\ell)\). The side equations hold to first order, and the derivative of the active diameter is
zero. Thus
\((v,\eta)\in\Kcal(p,\ell)\)
with
\(\eta>0\).
The direction is realised by
\(p_0(t)=p_0, \quad p_2(t)=p_2\),
\[
p_1(t)
=
\left(
\frac12,
-\frac25-t
\right),
\qquad
p_3(t)
=
\left(
\frac12,
\frac25+t
\right),
\]
and \(\ell(t)=\sqrt{1/4+(2/5+t)^2}\). 
For all sufficiently small \(t>0\), the active diameter remains equal to
\(1\), the other diagonal remains shorter than \(1\), and the polygon remains
strictly convex, so Proposition~\ref{prop:instability} applies.
Moreover, the multiplier identity \eqref{eq:multiplier-identity} cannot
hold at the rhombus, since neither \(p_1\) nor \(p_3\) is incident to an
active diameter and the incident side directions are linearly independent
at both vertices. The vertex balances therefore force
\(\alpha_0=\alpha_1=\alpha_2=\alpha_3=0\),
contradicting Remark~\ref{rem:multiplier-balance}.

\end{proof}

\section{The polygons \texorpdfstring{\(A_{16}\), \(A_{32}\) and \(A_{64}\)}{A16, A32 and A64}}
\label{sec:a16-data}

We recall the general \(A_n\) construction and record the geometric
data for \(A_{16}\), \(A_{32}\) and \(A_{64}\), giving explicit coordinates
and the diameter graph of \(A_{16}\). These data underpin the
certified first order analysis of the following section.

\subsection{The construction of the polygons \texorpdfstring{\(A_n\)}{An}}
\label{subsec:construction}

The OPTIGON file\footnote{\texttt{cstrt\_ngon/cstrt\_audet\_ngon.jl}}
implements the construction of Bingane and Audet~\cite{binganeaudet2022}
for \(n=2^k\), \(k\geq4\) \cite{optigon}.
For such \(n\), the constructor sets

\[
x(t)
=
\sin t
-
(2\cos t-1)
\frac{\cos((n/4-2)t)\sin(nt/4)}{\cos(2t)}.
\]

The constructor equation is
\(F_n(t)=\bigl(2x(t)+\sin((n/2-1)t)\bigr)^2+\cos^2((n/2-1)t)-4\sin^2(t/2)=0\),
and for \(n=16,32,64\) it has a unique root \(z_n\) in
\(I_{\mathrm{search}}=\bigl(2\arcsin(\sin(\pi/n)/2),\pi/n\bigr)\)
(Propositions~\ref{prop:a16-certified-checks}(1) and~\ref{prop:family}); the
interval is the root bracket of \cite{optigon}, while \cite{binganeaudet2022}
takes \(t\in(0,\pi/n)\). Then \(A_n\) has
side \(\ell=2\sin(z_n/2)\), vertices \(p_0=(0,0)\) and \(p_{n/2}=(0,1)\), and the
vertical line through them as axis of symmetry. Its diameter graph has \(n-1\) edges: the axis
\(\{0,n/2\}\), a path of length \(3n/8-1\) from \(0\) through \(n/2-1\) to
\(n/4\), the pendant edges \(\{4k-1,4k-1+n/2\}\) (\(k=1,\ldots,n/8\)), and the
images of the path and of the pendant edges under \(i\mapsto n-i\) \cite[Sect.~3]{binganeaudet2022}. We identify \(A_n\) with its pair \((p,\ell)\) and write, e.g., \(\Kcal(A_n)\) and \(\Ical(A_n)\)
for \(\Kcal(p,\ell)\) and \(\Ical(p,\ell)\).
Table~\ref{tab:polygons} and Figure~\ref{fig:a16} show \(A_{16}\), \(A_{32}\) and \(A_{64}\),
and Online Resource~1 lists the coordinates and the diameter graphs
of \(A_{32}\) and \(A_{64}\).

\begin{table}[ht]
\centering
\small
\caption{Constructor root and perimeter of the three polygons, rounded to the
digits shown.}
\label{tab:polygons}
\begin{tabular}{@{}rcc@{}}
\toprule
\(n\) & \(z_n\) & \(n\ell\)\\
\midrule
16 & \(0.1962703721\) & \(3.1352878881\)\\
32 & \(0.0981721090\) & \(3.1402460942\)\\
64 & \(0.0490872986\) & \(3.1412717079\)\\
\bottomrule
\end{tabular}
\end{table}

\subsection{The hexadecagon \texorpdfstring{\(A_{16}\)}{A16} in detail}
For \(n=16\), the construction is short enough to be written out in full, and
\(A_{16}\) carries the most detailed certificate of Section~\ref{sec:a16-certified},
including an explicit radius.

Here the numerator factor \(\cos((n/4-2)t)\) equals \(\cos(2t)\) and
cancels the denominator, provided \(\cos(2t)\neq0\); this holds on the root
interval used below, where \(\cos(2t)>0\). Since also \(nt/4=4t\), the
constructor reduces to
\(x(t) = \sin t-(2\cos t-1)\sin(4t)\). With
\(y(t)=\cos t-(2\cos t-1)\cos(4t)\), \(\bar x=x(t)\) and \(\bar y=y(t)\), it places the
vertices at
\begin{equation}
\label{eq:constructor}
\begin{aligned}
p_0&=(0,0), & p_1&=\bigl(\sin(2t)-\sin t,\ \cos t-\cos(2t)\bigr),\\
p_2&=\bigl(\sin(3t)-\sin t,\ \cos t-\cos(3t)\bigr), & p_3&=(-\bar x,\ \bar y),\\
p_4&=\bigl(\bar x+\sin(7t),\ \bar y+\cos(7t)\bigr), & p_5&=\bigl(\bar x+\sin(6t),\ \bar y+\cos(6t)\bigr),\\
p_6&=\bigl(\bar x+\sin(5t),\ \bar y+\cos(5t)\bigr), & p_7&=(\sin t,\ \cos t),\\
p_8&=(0,1), &&
\end{aligned}
\end{equation}
and the remaining vertices are given by the reflection \(p_{16-i}=(-x_i,y_i)\),
\(i=1,\ldots,7\), where \(p_i=(x_i,y_i)\). The constructor equation is \(F(t)=0\), where \(F=F_{16}\) is
\[
F(t)
=
\bigl(2x(t)+\sin(7t)\bigr)^2
+
\cos^2(7t)
-
4\sin^2(t/2)
=
|p_4-p_3|^2-|p_8-p_7|^2 ,
\]
so \(F(t)=0\) is the equality of the sides \(\{3,4\}\) and \(\{7,8\}\). It is the
only condition the construction leaves open, as verified in
Lemma~\ref{lem:exact-equalities}.

Evaluating the constructor formulas symbolically over \(\Q(c)[s]/(s^2+c^2-1)\), with \(c=\cos t\) and \(s=\sin t\) as indeterminates, reproduces the polynomial form of
Subsection~\ref{subsec:exact-algebra} and the constructor equation in the form
\(F_c(\cos t)\) given there \cite{chiqpkg}; the perimeter at the certified root agrees with
the value \(3.1352878881\) of \cite{binganeaudet2022} to ten decimal places. We write \(z^\ast=z_{16}\) for the
exact root, numerically \(z^\ast\approx0.1962703721303228\); certified enclosures of
\(\ell=2\sin(z^\ast/2)\) and \(16\ell\) are given in
Proposition~\ref{prop:a16-certified-checks}(2).

The polygon is symmetric about the vertical axis, with vertices \(0\)
and \(8\) on it; Table~\ref{tab:a16-coordinates} lists all \(16\) vertices,
and Figure~\ref{fig:a16}(a) shows the polygon. Decimal values in this section
are nominal and serve only to identify the polygon.
The declared active diameter graph is
\[
\begin{aligned}
E(A_{16})=\{&
\{0,7\},\{0,8\},\{0,9\},\{1,9\},\{2,9\},\\
&\{2,10\},\{3,10\},\{3,11\},\{3,12\},\{4,13\},\\
&\{5,13\},\{6,13\},\{6,14\},\{7,14\},\{7,15\}
\}.
\end{aligned}
\]

\begin{table}[ht]
\centering
\small
\caption{Nominal coordinates of \(A_{16}\), rounded to the digits shown.}
\label{tab:a16-coordinates}
\begin{tabular}{@{}r rr c r rr@{}}
\toprule
\(i\) & \(x_i\) & \(y_i\) & \hspace{1.5em} & \(i\) & \(x_i\) & \(y_i\)\\
\midrule
0 & 0.0000000000 & 0.0000000000 && 8 & 0.0000000000 & 1.0000000000\\
1 & 0.1875244690 & 0.0568606083 && 9 & $-$0.1950126739 & 0.9808007224\\
2 & 0.3603600643 & 0.1491991822 && 10 & $-$0.3465227585 & 0.8565298509\\
3 & 0.4847268698 & 0.3006305295 && 11 & $-$0.4389707792 & 0.6837527728\\
4 & 0.4959501447 & 0.4962643541 && 12 & $-$0.4959501447 & 0.4962643541\\
5 & 0.4389707792 & 0.6837527728 && 13 & $-$0.4847268698 & 0.3006305295\\
6 & 0.3465227585 & 0.8565298509 && 14 & $-$0.3603600643 & 0.1491991822\\
7 & 0.1950126739 & 0.9808007224 && 15 & $-$0.1875244690 & 0.0568606083\\
\bottomrule
\end{tabular}
\end{table}

\section{Certified verification of \texorpdfstring{\(A_{16}\), \(A_{32}\), \(A_{64}\), \(Z_{32}\) and \(Z_{64}\)}{A16, A32, A64, Z32 and Z64}}
\label{sec:a16-certified}

We verify the geometric and first order hypotheses of
Theorem~\ref{thm:first-order-topological} for the three
polygons of Section~\ref{sec:a16-data} and for the competing
configurations \(Z_{32}\) and \(Z_{64}\). We also obtain an
explicit neighbourhood of strict local maximality for \(A_{16}\).

\subsection{Exact identities and validated estimates}
\label{subsec:scope}

\begin{sloppypar}
The computation establishes, for the exact configuration \(A_{16}\), the
statements of Table~\ref{tab:certificate}: the two hypotheses of
Theorem~\ref{thm:first-order-topological} and the constants entering
Theorem~\ref{thm:a16-radius}. Equalities are proved in exact arithmetic over
\(\Q[c]\), where \(c=\cos t\); strict inequalities and regularity statements by
outward rounded interval arithmetic at \(256\)-bit precision, with the interval
Newton operator, a Krawczyk step \cite{krawczyk}, a preconditioned regularity
test and an interval Cholesky factorisation \cite{rump,alefeldmayer}. An interval
is never used to conclude an equality. The software is
\texttt{IntervalArithmetic.jl}~0.22.36 \cite{intervalarithmeticjl} under Julia~1.12.
For \(A_{16}\), \(A_{32}\), \(A_{64}\), \(Z_{32}\) and \(Z_{64}\) the linear algebra and every
evaluation over the certification boxes run in \(256\)-bit interval arithmetic;
in the general certifier of Subsection~\ref{subsec:family}, the isolation of
the root of \(A_{16}\), \(A_{32}\) and \(A_{64}\) uses
\(304\), \(368\) and \(496\) bits, respectively. Floating point arithmetic
only proposes starting data (the nominal root, the approximate inverses
used as preconditioners, the submatrix to test and the multiplier to
normalise); every conclusion is drawn from interval enclosures.

The accompanying Supplementary Material provides the detailed arithmetic
data, numerical bounds, geometric configurations and instructions for
reproducing the computations. The code and a Lean~4 \cite{leanfour}
formalisation are available in \cite{chiqpkg}. The formalisation covers
Theorem~\ref{thm:first-order-topological} and the results for
\(A_{16}\), \(A_{32}\) and \(A_{64}\)
(for instance \nolinkurl{ESP1.Family.A32.published_locally_maximal})
and, at every root of the defining equation in a certified bracket,
for \(Z_{32}\) and \(Z_{64}\), supplied as an additional check
on which no proof here depends. It rederives the hypotheses of
Proposition~\ref{prop:a16-certified-checks} with weaker constants and
a different positive multiplier vector, and it also proves
Theorem~\ref{thm:a16-radius}, once with constants of its own and once
by the argument of Subsection~\ref{subsec:radius} from the displayed
constants.

Subsection~\ref{subsec:family} extends the certificate to \(A_{32}\)
and \(A_{64}\), and Subsection~\ref{subsec:Z} to \(Z_{32}\) and
\(Z_{64}\). The general certifier of \cite{chiqpkg} runs the same
checks for all five polygons, and an independent implementation in
Python's standard library (own integer interval arithmetic, no Julia)
certifies the numerical hypotheses for all five and the exact
identities of \(Z_{32}\), \(Z_{64}\). The uniqueness of the parameter
\(t_n\) of \(Z_n\) (Subsection~\ref{subsec:Z}) rests on the interval
bisection of the Julia certificate; it is also proved in Lean
(for \(Z_{64}\), among the configurations that exist).

The package also runs the general certifier on \(A_n\) for
\(n\in\{128,256,512,1024\}\). We do not state those results,
because for \(n\geq65\) its rank and multiplier steps use floating
point interval products, a wider trust base than the \(256\)-bit
arithmetic used for the five polygons treated here.
\end{sloppypar}

\begin{table}[ht]
\centering
\small
\caption{The certificate for \(A_{16}\), with bounds rounded in the direction
that preserves them: \(B=D(e,d)\) (Definition~\ref{def:jacobian}), \(M\) the \(30\times30\) submatrix from the proof of Proposition~\ref{prop:a16-certified-checks}(7), \(M_0\) the
normalised system of Subsection~\ref{subsec:omitted-stationarity}, \(C\),
\(C_0\) approximate inverses of their midpoints, \(s_{30}(B)\) the least nonzero
singular value of \(B\), and \(V\) the matrix with columns \((\tau_i,0)\) from the
proof of Lemma~\ref{lem:euclidean-slice}; \(\Omega_{\mathrm{root}}\), \(\Omega_{\mathrm{cert}}\) are the boxes of
Subsection~\ref{subsec:boxes}, and \(N\) is the interval Newton operator of
Proposition~\ref{prop:a16-certified-checks}(1).}
\label{tab:certificate}
\begin{tabular}{@{}lll@{}}
\toprule
Statement & Certified value & Method\\
\midrule
unique root \(z^\ast\in \Omega_{\mathrm{root}}\) & \(N(\Omega_{\mathrm{root}})\subset\operatorname{int}\Omega_{\mathrm{root}}\) & interval Newton\\
no other root in \(I_{\mathrm{search}}\) & \(0\notin F'(I_{\mathrm{search}})\) & monotonicity\\
\(\Omega_{\mathrm{root}}\subset\operatorname{int}\Omega_{\mathrm{cert}}\) & verified & containment\\
\(\min\det(p_{i+1}-p_i,\,p_j-p_i)\) & \(>0.00748\) & \(224\) enclosures\\
\(\delta\) & \(>0.00809\) & \(105\) enclosures\\
\(\|\mathrm{Id}-C_0M_0\|_\infty\) & \(<2\times10^{-25}\) & regularity of \(M_0\)\\
\(\min\lambda_{ij}\), \(\|\alpha\|_2\) & \(>0.09186\), \(<4.004\) & Krawczyk step\\
\(\|\mathrm{Id}-CM\|_\infty\) & \(<1.5\times10^{-25}\) & regularity of \(M\)\\
\(s_{30}(B)^2\) & \(\geq0.09119\) & interval Cholesky of \(B^TB+VV^T-0.09119\,\mathrm{Id}\)\\
\bottomrule
\end{tabular}
\end{table}

\subsection{Exact algebra over \texorpdfstring{\(\Q[c]\)}{Q[c]}}
\label{subsec:exact-algebra}
Let \(T_j\) and \(U_j\) be the Chebyshev polynomials of the first and second kind. Since
\(\sin(jt)=\sin t\,U_{j-1}(\cos t)\) and \(\cos(jt)=T_j(\cos t)\), every vertex of \eqref{eq:constructor} has the form
\(p_k=\bigl(\sin t\,X_k(\cos t),Y_k(\cos t)\bigr)\) with
\(X_k,Y_k\in\Q[c]\), \(c\) an indeterminate. Hence
\(|p_i-p_j|^2=P_{ij}(\cos t)\) with
\(P_{ij}=(1-c^2)(X_i-X_j)^2+(Y_i-Y_j)^2\in\Q[c]\), and
\(F(t)=F_c(\cos t)\), where
\[
\begin{aligned}
F_c(c)={}&4096c^{12}-2048c^{11}-12288c^{10}+6656c^9
+12800c^8-7680c^7\\
&-4992c^6+3680c^5+192c^4-592c^3+192c^2-14c-1.
\end{aligned}
\]
\begin{lemma}
\label{lem:exact-equalities}
Let \(\Delta_i=P_{i,i+1}-(2-2c)\) for the side from \(p_i\) to \(p_{i+1}\), and
\(\Delta_{\{i,j\}}=P_{ij}-1\) for every declared active diameter edge \(\{i,j\}\). Then
\(\Delta_3=\Delta_{12}=F_c\) in \(\Q[c]\) (the sides \(\{3,4\}\) and \(\{12,13\}\)), and the remaining \(29\) polynomials are
identically zero. Consequently, if \(t\in I_{\mathrm{search}}\) is a root of \(F\), all \(16\) side equalities, with
\(\ell=2\sin(t/2)\), and all \(15\) declared active diameter equalities hold
exactly.
\end{lemma}
\begin{proof}
The identities follow by expansion in \(\Q[c]\) in exact integer arithmetic.
Since \(\ell^2=4\sin^2(t/2)=2-2\cos t\), the vanishing of \(\Delta_i(\cos t)\) is the
side equality, and \(F_c(\cos t)=F(t)=0\).
\end{proof}

\subsection{Root isolation and interval bounds}
\label{subsec:boxes}
The root is isolated in the interval \(\Omega_{\mathrm{root}}\) of values of \(t\), of radius \(10^{-30}\) about
a nominal root \(z\), a \(38\)-digit approximation of \(z^\ast\) given in the Supplementary Material.
All other quantities except \(\ell\) and \(16\ell\) are enclosed over the interval \(\Omega_{\mathrm{cert}}\), of radius
\(10^{-28}\) about the same \(z\). Since \(\Omega_{\mathrm{root}}\subset\operatorname{int}\Omega_{\mathrm{cert}}\),
every bound established over \(\Omega_{\mathrm{cert}}\) holds at the exact root
\(z^\ast\). These radii lie far below the resolution of double precision,
necessitating multiprecision interval arithmetic; every margin in
Table~\ref{tab:certificate} exceeds the resulting enclosure widths
by many orders of magnitude.

\subsection{Recovery of the omitted stationarity equations}
\label{subsec:omitted-stationarity}
Let
\[
\mu=\bigl((\alpha_i)_{i=0}^{15},(\lambda_{ij})_{\{i,j\}\in E(A_{16})}\bigr)
\]
be the multiplier vector, with \(\alpha=(\alpha_i)\) and \(\lambda=(\lambda_{ij})\), so that, identifying covectors with vectors by the inner product, \(DL\) becomes its
gradient \((0,\ldots,0,16)\) and
\eqref{eq:multiplier-identity} reads \(B^T\mu=DL\). Let \(M_0\) be the matrix of
the square system formed by the normalisation \(\lambda_{\{0,9\}}=\tfrac15\) and the \(30\)
stationarity equations for all ambient coordinates except \(x_1\), \(y_1\),
\(x_8\). When \(M_0\) is invertible, as established in
Proposition~\ref{prop:a16-certified-checks}(5), the system has
a unique solution \(\mu\), enclosed by a Krawczyk step (the full enclosure is in
\cite{chiqpkg}, only the bounds of Proposition~\ref{prop:a16-certified-checks}(6) are used in Subsection~\ref{subsec:radius}). The residual
\(\rho=B^T\mu-DL\in\R^{33}\) is then supported in the omitted coordinates:
\(\rho=\rho_{x_1}e_{x_1}+\rho_{y_1}e_{y_1}+\rho_{x_8}e_{x_8}\), with \(e_{x_1},e_{y_1},e_{x_8}\) the coordinate vectors.

Both \(B^T\mu\) and \(DL\) annihilate \(R_{\mathrm{eq}}(A_{16})\). If
\(\zeta\in R_{\mathrm{eq}}(A_{16})\), then \(B\zeta=0\), so
\(\langle B^T\mu,\zeta\rangle=\langle\mu,B\zeta\rangle=0\), and \(\zeta\) has zero side
length component, so \(DL[\zeta]=0\). Applying
\(\langle\rho,\zeta\rangle=0\) to the two translations and the rotation gives
\[
Q
\begin{pmatrix}
\rho_{x_1}\\
\rho_{y_1}\\
\rho_{x_8}
\end{pmatrix}
=
0,
\qquad
Q
=
\begin{pmatrix}
1&0&1\\
0&1&0\\
-y_1&x_1&-y_8
\end{pmatrix}.
\]
By \eqref{eq:constructor} and Lemma~\ref{lem:exact-equalities},
\(\det Q=y_1-y_8=y_1-1\leq|p_1-p_0|-1=\ell-1<0\), since \(\ell=2\sin(z^\ast/2)<1\). Hence \(\rho=0\) and
\(B^T\mu=DL\) holds exactly; beyond the root bracket of Proposition~\ref{prop:a16-certified-checks}(1), no interval arithmetic enters this step. Since \(B\) has
\(31\) rows, once its rank is shown to be \(30\) (Proposition~\ref{prop:a16-certified-checks}(7)), the multipliers
form a one parameter affine family, of which \(\mu\) is the representative with
\(\lambda_{\{0,9\}}=\tfrac15\).

\subsection{Certified geometric and first order data}

\begin{proposition}
\label{prop:a16-certified-checks}
The exact constructor root \(z^\ast\) and the configuration \(A_{16}\) satisfy:
\begin{enumerate}[label=\textup{(\arabic*)}]
\item \(F\) has a unique root \(z^\ast\in \Omega_{\mathrm{root}}\), and
\(\Omega_{\mathrm{root}}\subseteq I_{\mathrm{search}}\); moreover \(F'\) does not vanish on
\(I_{\mathrm{search}}\), so \(F\) is strictly monotone there and \(z^\ast\) is its
only root in \(I_{\mathrm{search}}\). This selects \(\cos z^\ast\) among the nine real roots of \(F_c\) in \([-1,1]\),
which are counted exactly, by a Sturm sequence, in the test suite of \cite{chiqpkg}.
\item At the exact root,
\[
\begin{gathered}
\ell\in[0.1959554930055385,\,0.1959554930055386],\\
16\ell\in[3.135287888088616,\,3.135287888088617].
\end{gathered}
\]
\item All \(224\) support determinants
\(\det\bigl(p_{i+1}-p_i,\;p_j-p_i\bigr)\), \(j\notin\{i,i+1\}\), exceed
\(0.00748\). Hence \(A_{16}\) is strictly convex, with the displayed labelled
cyclic order.
\item \(E(A_{16})\) is exactly as displayed, and
\(\delta(A_{16})>0.00809\).
\item The normalised system of Subsection~\ref{subsec:omitted-stationarity} is
regular: \(\|\mathrm{Id}-C_0M_0\|_\infty<2\times10^{-25}\); in particular \(M_0\) is
invertible. This item is established by interval evaluation over
\(\Omega_{\mathrm{cert}}\) alone, independently of
Subsection~\ref{subsec:omitted-stationarity}.
\item The identity \eqref{eq:multiplier-identity} holds with \(\lambda_{ij}>0.09186\)
for every \(\{i,j\}\in E(A_{16})\) and \(\|\alpha\|_2<4.004\); by Remark~\ref{rem:multiplier-balance}, \(\sum_{i=0}^{15}\alpha_i=-16\).
\item \(\rank B=30\); the \(30\times30\) submatrix \(M\) of the proof satisfies
\(\|\mathrm{Id}-CM\|_\infty<1.5\times10^{-25}\).
\item The least nonzero singular value of \(B=D(e,d)\) satisfies \(s_{30}(B)^2\geq0.09119\).
\end{enumerate}
\end{proposition}
\begin{proof}
(1) Interval Newton on \(\Omega_{\mathrm{root}}\): \(0\notin F'(\Omega_{\mathrm{root}})\)
and \(N(\Omega_{\mathrm{root}})=\check z-F(\check z)/F'(\Omega_{\mathrm{root}})\subset
\operatorname{int}\Omega_{\mathrm{root}}\), with \(\check z\in \Omega_{\mathrm{root}}\) its floating point midpoint; the endpoints
of \(\Omega_{\mathrm{root}}\) lie inside \(I_{\mathrm{search}}\); and an enclosure
of \(F'\) on \(I_{\mathrm{search}}\) excludes zero, so \(F\) is strictly
monotone there.
(2) Evaluation of \(2\sin(t/2)\) over \(\Omega_{\mathrm{root}}\). (3), (4)
Evaluation of the coordinate functions over \(\Omega_{\mathrm{cert}}\), together
with Lemma~\ref{lem:exact-equalities} for the \(15\) active pairs.

(5) The bound is the interval evaluation recorded in Table~\ref{tab:certificate};
since \(\|\mathrm{Id}-C_0M_0\|_\infty<1\), every matrix in the interval enclosure of
\(M_0\) is invertible, and so is \(M_0\). (6) Granted (5), the construction of
Subsection~\ref{subsec:omitted-stationarity} produces \(\mu\) and shows that
\(B^T\mu=DL\) holds exactly, which is \eqref{eq:multiplier-identity}; the bounds on
\(\lambda\) and \(\alpha\) are those of the Krawczyk enclosure of
Table~\ref{tab:certificate}, and \(\sum_{i=0}^{15}\alpha_i=-16\) is
Remark~\ref{rem:multiplier-balance}. (7) The submatrix \(M\) of \(B\) obtained by deleting the row of \(d_{\{0,9\}}\) and the columns \(x_1,y_1,x_8\) is
invertible, since \(\|\mathrm{Id}-CM\|_\infty<1\)
makes every matrix in its interval enclosure invertible, so \(\rank B\geq30\), and \(\rank B\leq2n-2=30\) by
Definition~\ref{def:jacobian}. (8) Let \(V\) be the matrix whose columns are the
vectors \((\tau_i,0)\) of the proof of Lemma~\ref{lem:euclidean-slice}; by (7) and
Definition~\ref{def:jacobian}, \(\operatorname{im}V=R_{\mathrm{eq}}(A_{16})=\ker B\).
The symmetric matrix \(B^TB+VV^T\) preserves the decomposition
\(R_{\mathrm{eq}}(A_{16})\oplus R_{\mathrm{eq}}(A_{16})^\perp\) and acts on the second
summand as \(B^TB\), so its least eigenvalue is at most \(s_{30}(B)^2\). The Cholesky
factorisation of \(B^TB+VV^T-0.09119\,\mathrm{Id}\), run in interval arithmetic on the
lower triangle of its enclosure, has every pivot enclosure positive; by inclusion
monotonicity the factorisation of every symmetric matrix in the enclosure runs to
completion with positive pivots, so each is positive definite \cite{alefeldmayer}.
Hence \(0.09119\) is smaller than the least eigenvalue of \(B^TB+VV^T\), and
therefore \(s_{30}(B)^2\geq0.09119\).
\end{proof}
\subsection{Strict local maximality of \texorpdfstring{\(A_{16}\)}{A16}}
\begin{theorem}
\label{thm:a16-certified}
The hexadecagon \(A_{16}\) is a strict local maximiser of
perimeter modulo orientation preserving Euclidean motions among convex
equilateral small hexadecagons.
\end{theorem}
\begin{proof}
By Proposition~\ref{prop:a16-certified-checks}(7), \(\rank B=30=2n-2\), so
\(\Ical(A_{16})=R_{\mathrm{eq}}(A_{16})\) by Definition~\ref{def:jacobian}.
Proposition~\ref{prop:a16-certified-checks}(6) gives
\eqref{eq:multiplier-identity} with every active diameter multiplier strictly
positive. By Lemma~\ref{lem:exact-equalities} and
Proposition~\ref{prop:a16-certified-checks}(3),(4), \(A_{16}\) is convex, equilateral and small.
Theorem~\ref{thm:first-order-topological} now applies.
\end{proof}
\subsection{An explicit neighbourhood}
\label{subsec:radius}
Because the squared constraints are quadratic, the certified data give a
neighbourhood of explicit size. We measure distance modulo motions by
\begin{equation}
\label{eq:dist}
\operatorname{dist}\bigl((q,m),(p,\ell)\bigr)
=
\min_{g}\Bigl(\sum_{i=0}^{15}|g(q_i)-p_i|^2+(m-\ell)^2\Bigr)^{1/2},
\end{equation}
the minimum over orientation preserving Euclidean motions \(g\). For each
rotation the best translation matches the centroids, so the minimum is taken
over the compact rotation group and is attained.
\begin{theorem}
\label{thm:a16-radius}
Let \((p,\ell)=A_{16}\), and let \((q,m)\in(\R^2)^{16}\times\R\) satisfy
\(|q_{i+1}-q_i|=m\) for every \(i\), \(|q_i-q_j|\leq1\) for all \(i,j\), and
\(m\geq\ell\); no convexity or cyclic order is assumed. If
\(\operatorname{dist}\bigl((q,m),(p,\ell)\bigr)<6.63\times10^{-4}\), then
\(m=\ell\) and \(q\) is the image of \(p\) under an orientation preserving
Euclidean motion.
\end{theorem}
\begin{proof}
Let \(g\) attain the minimum and write \(h=(\xi,\eta)=(gq-p,\,m-\ell)\); since
\(g\) preserves every distance, we may replace \(q\) by \(gq\). The first order
conditions of the minimisation are \(\sum_i\xi_i=0\) and
\(\sum_i\langle \xi_i,J(p_i+\xi_i)\rangle=0\). As \(\langle \xi_i,J\xi_i\rangle=0\), the
second is \(\sum_i\langle \xi_i,Jp_i\rangle=0\), so \(h\) is orthogonal to
\(R_{\mathrm{eq}}(A_{16})\).

Let \(B=D(e,d)\) be the Jacobian of Definition~\ref{def:jacobian} at \(A_{16}\),
with blocks \(B_e\) (sides) and \(B_d\) (active diameters). By
Lemma~\ref{lem:exact-equalities}, \(|p_{i+1}-p_i|=\ell\) and \(|p_j-p_i|=1\) on
the active pairs, and the squared constraints are quadratic, so their Taylor
expansions at \((p,\ell)\) are exact:
\[
|q_{i+1}-q_i|^2-m^2=2\ell\,(B_eh)_i+w_i,\qquad
|q_j-q_i|^2-1=2\,(B_dh)_{ij}+|\xi_j-\xi_i|^2,
\]
with \(w_i=|\xi_{i+1}-\xi_i|^2-\eta^2\). The left hand side of the first vanishes and that of the second is
nonpositive, so \(B_eh=-w/(2\ell)\) and every \((B_dh)_{ij}\leq0\). Moreover
\(\|w\|_2\leq\sum_i|\xi_{i+1}-\xi_i|^2+4\eta^2\leq4\|h\|^2\), since
\(\|(\eta^2,\ldots,\eta^2)\|_2=4\eta^2\) in \(\R^{16}\) and
\(|\xi_{i+1}-\xi_i|^2\leq2|\xi_{i+1}|^2+2|\xi_i|^2\).

By \eqref{eq:multiplier-identity}, \(16\eta=\alpha\cdot B_eh+\lambda\cdot B_dh\),
hence \(\lambda\cdot B_dh=16\eta+\alpha\cdot w/(2\ell)\geq-\|\alpha\|_2\|w\|_2/(2\ell)\)
since \(\eta\geq0\). All terms on the left are nonpositive, so, with
\(\lambda_{\min}=\min_{\{i,j\}}\lambda_{ij}\),
\(\lambda_{\min}\|B_dh\|_1\leq\|\alpha\|_2\|w\|_2/(2\ell)\) and
\[
\|Bh\|_2\leq\|B_eh\|_2+\|B_dh\|_1\leq\frac{2}{\ell}\Bigl(1+\frac{\|\alpha\|_2}{\lambda_{\min}}\Bigr)\|h\|^2 .
\]

On the other hand, \(\ker B=\Ical(A_{16})=R_{\mathrm{eq}}(A_{16})\)
by Proposition~\ref{prop:a16-certified-checks}(7) and Definition~\ref{def:jacobian}, and \(h\perp\ker B\), so
\(\|Bh\|_2\geq s_{30}(B)\,\|h\|\).
Therefore, for every \(\sigma>0\) with \(\sigma^2\leq s_{30}(B)^2\), either \(h=0\) or
\[
\|h\|\ \geq\ r:=\frac{\sigma\,\ell}{2\,(1+\|\alpha\|_2/\lambda_{\min})}.
\]
By Proposition~\ref{prop:a16-certified-checks}(8), \(\sigma^2=0.09119\) is admissible. Moreover,
\(\|\alpha\|_2<4.004\) and \(\lambda_{\min}>0.09186\) by
Proposition~\ref{prop:a16-certified-checks}(6), and
\(\ell>0.1959\) by Proposition~\ref{prop:a16-certified-checks}(2); hence
\(r\geq6.63\times10^{-4}\). Under the hypothesis of the theorem, \(\|h\|<r\), so \(h=0\).
\end{proof}

\subsection{The polygons \texorpdfstring{\(A_{32}\) and \(A_{64}\)}{A32 and A64}}
\label{subsec:family}
For \(n\in\{32,64\}\) the same steps apply to the polygons of
Section~\ref{sec:a16-data}, with the changes listed in the proof of
Proposition~\ref{prop:family}. Run in exact arithmetic, the
constructor again gives vertices \(\bigl(\sin t\,X_k(c),Y_k(c)\bigr)\) with
\(X_k,Y_k\in\Q[c]\), and, where \(\cos(2t)\neq0\), \(F_n(t)\) is a positive constant multiple of \(F_{c,n}(\cos t)\) for some
\(F_{c,n}\in\Q[c]\).

\begin{table}[ht]
\centering
\small
\caption{Certified margins for \(A_{32}\) and \(A_{64}\), rounded in the
direction that preserves them; the quantities are those of Table~\ref{tab:certificate},
with \(M\) and \(C\) the submatrix and approximate inverse of the proof of
Proposition~\ref{prop:family}, where the multipliers are pinned.}
\label{tab:family}
\begin{tabular}{@{}rllll@{}}
\toprule
\(n\) & \(\min\det\) & \(\delta\) & \(\min\lambda_{ij}\) & \(\|\mathrm{Id}-CM\|_\infty\)\\
\midrule
32 & \(>9.43\times10^{-4}\) & \(>9.62\times10^{-4}\) & \(>2.32\times10^{-2}\) & \(<2\times10^{-24}\)\\
64 & \(>1.18\times10^{-4}\) & \(>1.18\times10^{-4}\) & \(>5.36\times10^{-3}\) & \(<1\times10^{-23}\)\\
\bottomrule
\end{tabular}
\end{table}

\begin{proposition}
\label{prop:family}
Let \(n\in\{32,64\}\). The constructor equation has a unique root \(z_n\) in
\(I_{\mathrm{search}}\), and at
\(z_n\):
\begin{enumerate}[label=\textup{(\arabic*)},nosep]
\item all \(n\) side equalities and all \(n-1\) declared active diameter
equalities hold exactly; all but two of the \(2n-1\) squared residuals \(\Delta\) (as in
Lemma~\ref{lem:exact-equalities}) vanish
identically, and the other two vanish at every root of \(F_{c,n}\);
\item all \(n(n-2)\) support determinants are positive, so \(A_n\) is strictly
convex, with the labelled cyclic order of the constructor;
\item \(E(A_n)\) is exactly the edge set of the diameter graph of Subsection~\ref{subsec:construction};
\item the identity \eqref{eq:multiplier-identity} holds with every \(\lambda_{ij}>0\);
\item \(\rank B=2n-2\);
\item \(3.1402460942<32\ell<3.1402460943\) for \(n=32\) and
\(3.1412717079<64\ell<3.1412717080\) for \(n=64\), with \(\ell=2\sin(z_n/2)\).
\end{enumerate}
The certified margins are those of Table~\ref{tab:family}.
\end{proposition}
\begin{proof}
As for Proposition~\ref{prop:a16-certified-checks}, with these changes. The
root is isolated by interval Newton in a box \(\Omega_{\mathrm{root}}\) of radius \(10^{-30}\) inside
\(I_{\mathrm{search}}\), and \(0\notin F_{c,n}'(\cos I_{\mathrm{search}})\). For
(5), a \((2n-2)\times(2n-2)\) submatrix \(M\) of \(B\), with its approximate inverse \(C\), is certified regular
(Table~\ref{tab:family}); (6) is the evaluation of \(2n\sin(t/2)\) over \(\Omega_{\mathrm{root}}\). The rank and multiplier steps run in \(256\)-bit interval arithmetic throughout; the isolation of the root uses
\(368\) bits for \(n=32\) and \(496\) for \(n=64\). For (4), the stationarity equations are kept except
at the gauge coordinates \(x_0,y_0,x_j\) (\(j=3\) for \(n=32\), \(j=5\) for
\(n=64\)), and the multiplier of largest modulus in a floating point least squares solution is fixed at its floating point value; the square system is certified
regular and solved by a Krawczyk step, and \(Q\), with rows \((1,0,1)\), \((0,1,0)\), \((-y_0,x_0,-y_j)\) and \(\det Q=y_0-y_j=-y_j\), is certified regular by the preconditioned test, which recovers
the three omitted equations.
\end{proof}

\begin{theorem}
\label{thm:family}
For \(n\in\{32,64\}\), the polygon \(A_n\) is a strict local
maximiser of perimeter modulo orientation preserving Euclidean motions among
convex equilateral small \(n\)-gons.
\end{theorem}
\begin{proof}
By Proposition~\ref{prop:family}(5) and Definition~\ref{def:jacobian},
\(\Ical(A_n)=R_{\mathrm{eq}}(A_n)\). Proposition~\ref{prop:family}(4)
gives positive multipliers, and by Proposition~\ref{prop:family}(1)--(3), \(A_n\) is convex,
equilateral and small, so Theorem~\ref{thm:first-order-topological}
applies.
\end{proof}

\subsection{The polygons \texorpdfstring{\(Z_{32}\) and \(Z_{64}\)}{Z32 and Z64}}
\label{subsec:Z}
The polygons \(Z_n\) of the introduction are described in
\cite[Props.~1 and~2]{binganeaudet2022}, which state that they are small, convex and of larger perimeter
than \(A_n\). Bingane and Audet publish no code for them, so we rebuild them from those
statements; Online Resource~1 lists the vertices and the diameter graphs. Write \(S_k=\sin kt\), \(C_k=\cos kt\). For
\(Z_{32}\), the diameter graph has the axis \(\{0,16\}\), the half path
\(0\)--\(11\)--\(24\)--\(10\)--\(23\), the pendant edges \(\{0,k\}\)
(\(12\leq k\leq15\)) and \(\{11,k\}\) (\(25\leq k\leq31\)), and their images under
\(i\mapsto n-i\); the parameter \(t\) satisfies
\[
G_{32}(t)=\bigl(2(S_5-S_{13}+S_{14})-S_{15}\bigr)^2+C_{15}^2-(2-2\cos t)=0 .
\]
For \(Z_{64}\) the half path visits, in order, the vertices
\[
0,31,63,30,61,29,60,28,58,27,57,26,56,25,55,24,54,23,53,21,52,19,51,16,
\]
and adds, at its successive edges, alternately \(+(S_m,C_m)\) and
\(-(S_m,C_m)\), beginning with \(+\), for \(m=1,2,3,\allowbreak5,6,7,8,\allowbreak10,\allowbreak\ldots,20,22,\allowbreak23,25,26\), so that it reaches \(p_{51}\), which depends on
\(t\) only; one more step, with \(m=29\), leads to \(p_{16}\), so that \(p_{16}=p_{51}+(S_{29},C_{29})\). The four vertices that carry the unknown \(y\) are
\(p_{50}=(-\tfrac12,y)\), \(p_{14}=(\tfrac12,y)\), \(p_{15}=(-\tfrac12+C_1,\,y+S_1)\) and \(p_{49}=(\tfrac12-C_1,\,y+S_1)\); two equations, \(|p_{51}-p_{50}|=|p_{16}-p_{15}|=2\sin(t/2)\),
fix \(t\) and \(y\). Write \(a=x_{51}+\tfrac12\) and let \(\ell^2=2-2\cos t\) be the squared side
length; the square root in what follows is real where
\(\ell^2\geq a^2\). The first equation gives \(y=y_{51}+\sqrt{\ell^2-a^2}\) (the
convex branch), and the second becomes
\[
G_{64}(t)=(a+S_{29}-C_1)^2+\bigl(C_{29}-S_1-\sqrt{\ell^2-a^2}\bigr)^2-\ell^2=0 .
\]

The diameter graph of \(Z_{64}\) consists of the \(23\) edges of the half path, including \(\{51,16\}\), and their images under \(i\mapsto n-i\),
the axis \(\{0,32\}\), the path \(15\)--\(50\)--\(14\)--\(49\) (\(3\) edges), and the pendant edges \(\{30,62\}\), \(\{28,59\}\), \(\{53,22\}\),
\(\{52,20\}\), \(\{51,18\}\), \(\{51,17\}\) (with their images under \(i\mapsto n-i\)). The paper of Bingane and Audet
gives no angle for the pendant edges, and none is needed. Each pendant vertex is
fixed, up to the reflection that convexity excludes, by the lengths its sides impose
(\(\ell\) to both neighbours for \(p_{62},p_{59},p_{22},p_{20}\), the diameter then following);
for \(p_{17}\) and \(p_{18}\), by \(\ell\) to the neighbour on the half path (\(p_{16}\), resp.\ \(p_{19}\)) and
\(1\) to \(p_{51}\). That the points \(p_{30}-(S_4,C_4)\),
\(p_{28}-(S_9,C_9)\), \(p_{53}+(S_{21},C_{21})\), \(p_{52}+(S_{24},C_{24})\), \(p_{51}+(S_{27},C_{27})\)
and \(p_{51}+(S_{28},C_{28})\) satisfy these conditions is part of the exact
verification of Proposition~\ref{prop:zfamily}(1).

\begin{table}[ht]
\centering
\small
\caption{Certified margins for \(Z_{32}\) and \(Z_{64}\), rounded in the direction
that preserves them, with columns as in Table~\ref{tab:family}; the multipliers of \(Z_{32}\) are pinned as in the proof of
Proposition~\ref{prop:family} (those of \(Z_{64}\) are unique). The values are
those of the certificate of Proposition~\ref{prop:zfamily}, and \(n\ell\) agrees with the
perimeter reported in \cite{binganeaudet2022}.}
\label{tab:zfamily}
\begin{tabular}{@{}rlllll@{}}
\toprule
\(n\) & \(\min\det\) & \(\delta\) & \(\min\lambda_{ij}\) & \(\|\mathrm{Id}-CM\|_\infty\) & \(n\ell\)\\
\midrule
32 & \(>9.43\times10^{-4}\) & \(>1.12\times10^{-4}\) & \(>4.06\times10^{-3}\) & \(<2\times10^{-24}\) & \(3.1403202339\ldots\)\\
64 & \(>1.18\times10^{-4}\) & \(>2.09\times10^{-5}\) & \(>2.80\times10^{-3}\) & \(<4\times10^{-22}\) & \(3.1412752155\ldots\)\\
\bottomrule
\end{tabular}
\end{table}

\begin{proposition}
\label{prop:zfamily}
Let \(n\in\{32,64\}\). The equation \(G_n(t)=0\) has exactly one root \(t_n\) in
\((0,\pi/n]\) at which \(Z_n\) exists, namely \(t_{32}=0.0981744286\ldots\) and
\(t_{64}=0.0490873533\ldots\), the values reported in \cite{binganeaudet2022}. For the
polygon \(Z_n\) built at \(t_n\):
\begin{enumerate}[label=\textup{(\arabic*)},nosep]
\item all \(n\) side equalities and all declared active diameter equalities hold
exactly: the squared residuals \(\Delta\) (as in Lemma~\ref{lem:exact-equalities}; for \(Z_{64}\),
polynomials in \(c\), \(s=\sin t\) and \(y\) modulo \(s^2+c^2-1\)) of the \(n-1\) (for
\(Z_{32}\)) or \(62\) (for \(Z_{64}\)) diameters vanish identically, those of all but two
(resp.\ four) sides vanish identically, and the others are the defining
equation(s);
\item all \(n(n-2)\) support determinants are positive, so \(Z_n\) is strictly
convex;
\item every other pair of vertices is at distance less than \(1\), so \(Z_n\) is
small and \(E(Z_n)\) is the declared edge set;
\item the identity \eqref{eq:multiplier-identity} holds with every \(\lambda_{ij}>0\);
\item \(\rank B=2n-2\);
\item \(3.1403202339<32\ell<3.1403202340\) for \(n=32\) and
\(3.1412752155<64\ell<3.1412752156\) for \(n=64\).
\end{enumerate}
The certified margins are those of Table~\ref{tab:zfamily}.
\end{proposition}
\begin{proof}
The proof follows that of Proposition~\ref{prop:family}, with these changes. For \(Z_{64}\), the four side
residuals that do not vanish are the two equations that define \(y\) and \(t\), and
they vanish at the polygon because \(y-y_{51}=\sqrt{\ell^2-a^2}\) and \(t_n\) is a root of \(G_{64}\). The
root is isolated by interval Newton in a box of radius \(10^{-30}\) about the nominal
value, with the closed form of \(G_n'\), which gives exactly one root in the box. Interval
bisection on \((0,\pi/n]\) discards every other box, because \(G_n\) is nonzero on that box or
the polygon does not exist (\(\ell^2<a^2\) for \(Z_{64}\) on most of the interval), so
\(t_n\) is the only root there. \(Z_{32}\) has \(2n-1\) active constraints of rank \(2n-2\), so its
multipliers form an affine family of dimension one and one of them is pinned, as
for \(A_n\); \(Z_{64}\) has \(2n-2\), and its multipliers are unique. In both
cases a square system is certified regular as in
Proposition~\ref{prop:family}. The certificate is that of \cite{chiqpkg}, with the
same checks as for \(A_{32}\) and \(A_{64}\) except that there is no published
constructor to run and uniqueness comes from the bisection. The independent implementation described in Subsection~\ref{subsec:scope} agrees with the Julia certificate provided in the accompanying Supplementary Material.
\end{proof}

\begin{theorem}
\label{thm:zfamily}
For \(n\in\{32,64\}\), the polygon \(Z_n\) is a convex equilateral small
\(n\)-gon and a strict local maximiser of perimeter modulo orientation preserving
Euclidean motions among convex equilateral small \(n\)-gons. Its perimeter exceeds
that of \(A_n\) by more than \(7.4\times10^{-5}\) for \(n=32\) and
\(3.5\times10^{-6}\) for \(n=64\); in particular \(A_{32}\) and \(A_{64}\) are not
global maximisers.
\end{theorem}
\begin{proof}
By Proposition~\ref{prop:zfamily}(5) and Definition~\ref{def:jacobian},
\(\Ical(Z_n)=R_{\mathrm{eq}}(Z_n)\), and
Proposition~\ref{prop:zfamily}(4) gives positive multipliers, so
Theorem~\ref{thm:first-order-topological} applies. By
Proposition~\ref{prop:zfamily}(1)--(3), \(Z_n\) is convex, equilateral and small.
The perimeters are those of Propositions~\ref{prop:zfamily}(6) and~\ref{prop:family}(6).
\end{proof}

\begin{remark}
The Lean development covers \(Z_{32}\) and \(Z_{64}\) as it does the other three polygons: at every
root of the defining equation in the certified bracket the polygon is strictly convex, equilateral,
small and a strict local maximiser, and a root exists there. It also proves that this root is the only one
in \((0,\pi/n]\). For \(Z_{32}\) the derivative of the defining polynomial in \(\cos t\) keeps one sign on
\([\cos(\pi/32),1]\). For \(Z_{64}\) the statement is about the configurations for which the construction
is defined, those with \(\ell^2\geq a^2\). There the proof uses that \(a=x_{51}+\tfrac12\) is positive and increasing in \(\cos t\), so that the
radicand decreases, and bounds the derivative of the defining equation in \(\cos t\) on subbrackets, the term coming from the
square root having the favourable sign; the equation is therefore strictly increasing on those configurations. It does not cover the identification of \(Z_n\)
with the polygon described in \cite{binganeaudet2022}, which rests on
its Propositions~1 and~2 as read above: the equation for \(t\), the half path, the
reflection, the diameter graph and the perimeter agree with them.
\end{remark}

\section*{Acknowledgements}
During the preparation of this article, we used OpenAI's
ChatGPT and Anthropic's Claude to explore mathematical arguments
guided by questions and conjectures formulated by the authors.
The assistance included the following activities:

\begin{itemize}

\item Exploring proof strategies for the strict local maximality
criterion (Theorem~\ref{thm:first-order-topological}),
including the treatment of dependent constraints and
infinitesimal Euclidean motions.

\item Assisting with the development and examination of
computational verification procedures, including exact algebra,
validated interval arithmetic, multiplier positivity, rigidity
checks, and the accompanying Lean~4 formalisation.

\item Exploring the derivation and verification of the explicit
neighbourhood of strict local maximality for $A_{16}$
(Theorem~\ref{thm:a16-radius}), including the use of
multiplier bounds and singular value estimates.

\item Identifying potential gaps in mathematical arguments,
suggesting corrections, and assisting with literature searches,
LaTeX preparation, and editorial revisions.

\end{itemize}

The authors checked the mathematical arguments and computational
results incorporated into the final manuscript and take full
responsibility for its content.

\end{document}